\documentclass{article}
\usepackage[left=2.4cm, right=2.5cm, top=2.5cm, bottom=2.5cm]{geometry}
\usepackage{graphicx} 
\usepackage{amsmath}
\usepackage{amssymb}
\usepackage{amsthm}
\usepackage{lineno}
\usepackage[stable]{footmisc}
\usepackage{blkarray}
\usepackage{enumerate}
\usepackage[shortlabels]{enumitem}
\usepackage[font=small,skip=0pt]{caption}
\usepackage{verbatim}
\usepackage{geometry}
\usepackage{algorithmicx}
\usepackage{float}
\usepackage{blindtext}
\usepackage{calc}
\usepackage{upgreek}
\usepackage{url}
\usepackage{setspace}

\usepackage{xcolor}
\definecolor{Maroon}{HTML}{ad2231}
\definecolor{webgreen}{HTML}{008000}
\usepackage{dsfont}
\usepackage{hyperref}
\usepackage{lipsum}
\hypersetup{colorlinks, breaklinks, urlcolor=black, linkcolor=black, citecolor=webgreen} 

\DeclareMathAlphabet{\mathpzc}{OT1}{pzc}{m}{it}

\numberwithin{equation}{section}

\newtheorem*{Note*}{Note}

\newtheorem*{Recall*}{Recall}

\usepackage{amsthm}
\newtheorem{theorem}{Theorem}

\newtheorem{proposition}[theorem]{Proposition}
\newtheorem{lemma}[theorem]{Lemma}
\newtheorem{remark}[theorem]{Remark}

\theoremstyle{definition}

\allowdisplaybreaks

\title{Fluctuations of  Omega-killed level-dependent spectrally negative L\'evy processes}
\author{
	Zbigniew Palmowski\footnote{Department of Applied Mathematics,
		Wroc\l{}aw Un. of Science and Technology, \texttt{zbigniew.palmowski@gmail.com}}\,, \,   Meral \c{S}im\c{s}ek\footnote{Institute of Applied Mathematics  \& Dept of
		Statistics, Middle East Technical University, \texttt{smeral@metu.edu.tr}} \; \&   Apostolos D. Papaioannou\footnote{Department of Mathematical Sciences, IFAM,
		University of Liverpool, \texttt{papaion@liverpool.ac.uk}}    
}

\newcommand{\email}[1]{\gdef\@email{\url{#1}}}

\newcommand{\N}{\mathbb{N}}

\newcommand{\R}{\mathbb{R}}

\newcommand{\Prob}{\mathbb{P}}
\newcommand{\Ex}{\mathbb{E}}

\DeclareMathAlphabet{\mathpzc}{OT1}{pzc}{m}{it}
\newcommand{\dd}{\mathrm}

\begin{document}

\maketitle

\tableofcontents

\vspace{0.1in}

\begin{abstract}
 In this paper, we solve exit problems  for a level-dependent L\'evy process which   is exponentially killed with a killing intensity that depends on the present state of the process. Moreover, we analyse the respective resolvents. All identities are given in terms of new generalisations of scale functions (counterparts of the scale function from the theory of L\'evy processes), which  are solutions of Volterra integral equations.  Furthermore, we obtain similar results for the reflected level-dependent L\'evy processes. The existence of the   solution of the  stochastic differential equation for reflected level-dependent L\'evy processes is also discussed. Finally, to illustrate our result, the probability of bankruptcy is obtained for an insurance risk process.
\end{abstract}

\noindent {\sc Keywords}: Level-dependent L\'evy processes; Omega model; Fluctuation theory; Volterra equation; Potential measures


\section{Introduction}
Exit problems for spectrally negative L\'evy processes (SNLP) have been extensively studied over the last 40 years and applied in many fields, such as financial mathematics, risk and queueing theory, physics and others. An overview can be found in \cite{Bertoin1996}, \cite{Kyprianou2014}. In more recent years, a parallel theory for fluctuations has been developed for more general L\'evy processes, such as  the refracted L\'evy processes introduced in \cite{kyprianou-refract}, or its generalisation in the frame of  the level-dependent L\'evy process introduced in \cite{CPRY2019}.
A level-dependent L\'evy process $U = (U_t)_{t\geq 0}$ is solution of the following SDE
	\begin{equation}
		\mathrm{d}U_t = \mathrm{d}X_t - \phi(U_t)\mathrm{d}t,
		\label{eqn:levelDependent}
	\end{equation}
where $X_t$ is a spectrally negative L\'evy process and $\phi$ is non-decreasing and locally Lipschitz continuous. The theory of level-dependent L\'evy processes was first introduced when $X$ is  a bounded variation L\'evy process in the setting of storage processes, see \cite{resnickJAP}, or risk processes, see Chapter VII.1 in \cite{Asmussen_Albrecher-2010}. In the latter case, the alternative form of Eq. \eqref{eqn:levelDependent}  can be
written as
	\begin{equation*}
		\mathrm{d}U_t =  -\mathrm{d}S_t + p(U_t)\mathrm{d}t,
	\end{equation*}
where $S_t$ is a compound Poisson process and
 $p(U_t) = ct - \phi(U_t)$ is a level-dependent premium rate, which has numerous application in risk models with dividend payments.
A comprehensive  work regarding  the fluctuation  theory   of level-dependent process  in Eq. \eqref{eqn:levelDependent} can be found in \cite{CPRY2019}, in which the authors show that the solution  for exit problems and the   resolvents are expressed in terms of scale functions that satisfy Volterra integral equations of second kind.

The aim of this paper is to develop the theory of scale functions in level-dependent SNLP which is killed exponentially with intensity that depends on the level of the process  and generalise  known exit identities for $\omega$-killed (reflected) SNLP, for which previous results in \cite{CPRY2019}, \cite{boli2018} (and \cite{PY2018}) are special cases. To formulate mathematically our problem, let $\omega: \mathbb R \to \mathbb R^+$ be a non-negative, locally bounded  measurable function and
	\begin{equation*}
		\kappa_x^{+}:=\inf\{t>0: U_t \geq x\} \quad \text { and } \quad \kappa_x^{-}:=\inf\{t>0: U_t < x\}.
		\label{eqn:passageLevelD}
	\end{equation*}
Throughout the paper, we denote the law of $X$, such that $X_0=x$, by $\mathbb P_x$ and the corresponding expectation by $\mathbb E_x$, with the convention $\mathbb P$ and $\mathbb E$ when $x=0$. Our main interest is to derive explicit expressions  for the occupation times weight by $\omega$ of the form, for $x\in[0,a]$,
	\begin{equation*}
	\mathpzc{A}(x, a) := \Ex_x\Bigl[ e^{-\int_{0}^{\kappa_a^+}\omega(U_t)\mathrm{d}t}\mathbf{1}_{(\kappa_a^+< \kappa_0^-)}\Bigr], \quad  	\mathpzc{B}(x, a): = \Ex_x\Bigl[ e^{-\int_{0}^{\kappa_0^-}\omega(U_t)\mathrm{d}t}\mathbf{1}_{( \kappa_0^- <\kappa_a^+)}\Bigr],
\end{equation*}
as well as similar closed formulas for the reflected at infima level-dependent process and finally explicit identities  for the respective resolvents (both in the non-reflected and reflected case). Apart from the above, we point out that we provide a comprehensive study for the  reflected  level-dependent L\'evy process, which  has not studied in the in the literature before, which is achieved by the strong approximation.

The $\omega$-killed exit problems have been previously considered for SNLP in \cite{boli2018}, \cite{BoliXiaowen2018},  for diffusion processes in \cite{BL2000}, for Markov additive processes in \cite{CKP2020} and have  attracted attention as they have numerous applications in insurance and finance. For instance in finance, $\omega(\cdot)$ can model a sophisticated structure that captures the interest rates which depend on the economic environment changes (see \cite{BL2000}), whilst in insurance $\omega(\cdot)$ models the rate of bankruptcy in a risk model that allows recovery when surplus is negative, namely the Omega risk model (see  \cite{GSH2012}).

 In  this paper, in order to handle the complexity of the process in  \eqref{eqn:levelDependent} under the $\omega$-killing,   we extend and  amend   methodologies used in \cite{CPRY2019} \cite{boli2018}, \cite{BoliXiaowen2018}, and \cite{PY2018}. Nevertheless, new techniques and arguments, which shorten and simplify the proofs of the results are derived in a number of results (see  the establishment of the scale functions of Section  \ref{killed scale} through  Proposition \ref{pro2w},  Theorem \ref{thm:twoSidedBelowLevelw}, to mention a few). This also includes a new  strong approximation in Section \ref{reflected} (see Proposition \ref{prop:convergenceVn}) for the reflected level-dependent process,  achieved by Proposition  \ref{prop:ConvergenceVk}.

The paper is organised as follows. In Section \ref{prelim} we provide some preliminaries. In Section \ref{killed scale}, we provide the the $\omega$-killed scale function as solutions to integral equations. In the same section identities for the two-sided exit upwards/downwards problem and the potential measure are derived. To illustrate the applicability of the results  the probability of bankrupcty for an Omega level-dependent model is derived, which  generalises   the Omega-model for a  L\'evy risk process analysed in \cite{AGS2011, GSH2012, LRZ2011, LRZ2014}.  In Section \ref{reflected}, we introduce the reflected level-dependent L\'evy process and show the existence of  its solution (using associated tools  form the Appendix).   In the same section, the identities for the  the upward exit problem and the potential measure are also derived.
\section{Preliminaries}
\label{prelim}
Let $X = (X_{t})_{t \geq 0}$ be a SNLP defined on the filtered space $(\Omega, \mathcal{F}, (\mathcal{F}_{t})_{t \geq 0}, \mathbb{P})$, where the filtration $(\mathcal{F}_{t})_{t \geq 0}$ is assumed to satisfy the usual assumptions of right continuity and completion. A L\'evy process with no positive jumps has Laplace transform of the form
	\begin{equation*}
		\Ex[e^{\theta X_t}]=e^{t\psi(\theta)}, \;\;  \forall \theta\geqslant 0.
	\end{equation*}
where  $\psi(\theta):[0,\infty)\rightarrow \mathbb R$,  known  as the  Laplace exponent of $X$, is a  continuous and strictly convex function given by the L\'evy-Khintchine formula
	\[ \psi(\theta) = \mu\theta + \frac{\theta^2 \sigma^2}{2} +
	\int_{(-\infty, 0)}^{} (e^{\theta x} -1- \theta x \mathbf{1}_{\{x >- 1\}}) \nu(\mathrm{d}x),\]
where $\mu \in \R$, $\sigma \geq 0$
and $\nu$ is a measure
concentrated on $(-\infty,0)$ called the L\'evy measure that satisfies $\int_{(-\infty,0)}(1 \wedge |x|^2) \nu(\mathrm{d}x)  <\infty$ . It is well known that $X$ has paths of bounded variation if and only if $\sigma =0$, and $\int_{(-1,0)}\nu(\mathrm  dx)<\infty$. In this case $X$ can be written as
	$ X_t=ct-S_t$,
where $ c:=\mu -\int_{(-1,0)}x\nu(\mathrm d x)$ and $(S_t)_{t\geqslant 0}$ is a driftless subordinator. For more details the reader is refereed to \cite{Bertoin1996} or \cite{Kyprianou2014}.

It is well-known, see Chapter 8 in \cite{Kyprianou2014}, \cite{KK2012}, that fluctuation identities for SNLP rely on the so-called scale functions. For  any $q\geq 0$, the  two families of $q$-scale functions $W^{(q)}: \R \to [0, \infty)$ and $Z^{(q)}: \R \to [1, \infty)$  are defined by
		\begin{equation*}
			\int_{0}^{\infty} e^{-\theta x}W^{(q)}(x)\mathrm{d}x = \frac{1}{\psi(\theta) -q}, \;  \theta > \Phi(q), \quad 	Z^{(q)}(x) = 1 + q\int_{0}^{x}W^{(q)}(y)\mathrm{d}y,
			\label{eqn:scaleW}
		\end{equation*}
where  $W^{(q)}$ is  right continuous and increasing function with $W^q(x)=0$ for $x<0$ and $Z^{(q)}$ to  inherit the  properties of $W^{(q)}$. In the rest of the paper, when $q=0$, we shall drop the subscript.

The  theory   of the level-dependent L\'evy process $U$, in \eqref{eqn:levelDependent}, heavily depends on the assumption
\begin{itemize}\item[]{\textbf{[A]}}
	The non-decreasing function $\phi$ is either locally Lipschitz continuous or in the (so-called \textit{multi-refracted}) form $\phi_k(x) = \sum_{j=1}^{k}\delta_j \mathbf{1}_{(x>b_j)}$, with  $0 < \delta_1, \ldots, \delta_k$, and $-\infty< b_1 < \cdots < b_k< \infty$. Moreover, $\phi(x) = 0$ for $x \leq d'$, where $d' \in \R$ is fixed and
	for the bounded variation case, it is assumed that $\phi(x) < c$.
	\label{assumption}
\end{itemize}
Under the Assumption \textbf{[A]}, it can be shown that  for multi-refracted case ($\phi\equiv \phi_k$) and when $U$ is of bounded variation, a strong   solution of  SDE in Eq.  \eqref{eqn:levelDependent}  exists and can be extended by a strong approximation argument to the case that $U$ is of unbounded variation [see Theorem 1 in \cite{CPRY2019}].   This is the base for proving (using a limiting argument) the existence  for the  solution of the SDE in Eq. \eqref{eqn:levelDependent} for a general rate function $\phi$, as it given in the proposition below [see Proposition 20 in \cite{CPRY2019}].
\begin{proposition}
	\label{thm:ExistenceU}
	Suppose that the locally Lipschitz continuous rate function $\phi$ satisfies condition \textbf{[A]}. Then, there exists a  solution $U$ to the SDE \eqref{eqn:levelDependent} with rate function $\phi$.
\end{proposition}
\noindent Exit identities and the potential measures for the  level-dependent L\'evy processes $U$ are  established
in terms of level-dependent scale functions $\mathbb{W}^{(q)}(x;y)$ and $\mathbb{Z}^{(q)}(x;y)$, $x$, $y \in \R$ [see \cite{CPRY2019}],   which are  solutions to the following integral equations, provided that
	$\mathbb{W}^{(q)}$ and $\mathbb{Z}^{(q)}$ are almost everywhere differentiable
		\begin{align}
		\mathbb{W}^{(q)}(x; y) &= W^{(q)}(x-y)+\int_y^x W^{(q)}(x-z)
		\phi(z)\mathbb{W}^{(q)\prime}(z; y)\mathrm{d}z,
		\label{eqn:shiftedLevelScaleW}\\
		\mathbb{Z}^{(q)}(x;y) &= Z^{(q)}(x-y)+\int_y^x W^{(q)}(x-z)
		\phi(z)\mathbb{Z}^{(q) \prime}(z;y)\mathrm{d}z,
		\label{eqn:shiftedLevelScaleZ}
	\end{align}
with boundary conditions $\mathbb{W}^{(q)}(y;y) =W^{(q)}(0)$ and
$\mathbb{Z}^{(q)}(y;y) = 1$ where $W^{(q)}$  and  $Z^{(q)}$ as defined above. We denote
$\mathbb{W}^{(q)}(x) = \mathbb{W}^{(q)}(x; 0)$ and $\mathbb{Z}^{(q)}(x) = \mathbb{Z}^{(q)}(x; 0)$ (we point out that the special case $\mathbb{Z}^{(q)}(x)$ is used in \cite{CPRY2019}).
We note that the general derivative of $W^{(q)}$ is not defined for all $x \in \R$.
Hence,  $W^{(q)\prime}$, $\mathbb{W}^{(q)\prime}$ are  understood as the right-derivatives, which always exists
[see Lemma 2.3 in~\cite{KK2012}].
Then, for any $x \geq 0$,
$\mathbb{W}^{(q)\prime}$ and $\mathbb{Z}^{(q)\prime}$ are solutions to following the Voltera-type
equations
	\begin{align*}
		\mathbb{W}^{(q)\prime}(x;y)
		&=\Xi_\phi(x)^{-1} W^{(q)\prime}((x-y)+)+ \Xi_\phi(x)^{-1}\int_y^x  W^{(q)\prime}(x-z)\phi(z) \mathbb{W}^{(q)\prime}(z;y)\mathrm{d}y,
		\\
		\mathbb{Z}^{(q)\prime}(x;y)
		&=\Xi_\phi(x)^{-1} q W^{(q)}(x-y)+ \Xi_\phi(x)^{-1} \int_y^x  W^{(q)\prime}(x-z)\phi(z)
		\mathbb{Z}^{(q)\prime}(z;y) \mathrm{d}z,
	\end{align*}
with  $\Xi_\phi(x)$ be  strictly positive function by assumption \textbf{[A]}
and  defined as
	\begin{equation}
		\Xi_\phi(x):=1-W^{(q)}(0) \phi(x), \;\;  \text{with} \; \; 	W^{(q)}(0)=
		\begin{cases}
			0 & \text{if $X$ has unbounded variation,} \\
			1/c & \text{if $X$ has bounded variation.}
		\end{cases}
		\label{eqn:Xi}
	\end{equation}
We note that the above Volterra integral equations for $\mathbb{W}^{(q)\prime}$ and $\mathbb{Z}^{(q)\prime}$  have a unique solutions in terms of the the scale function $W^{(q)}$ of the driving L\'evy process $X$ (for more details see Proposition 26 in \cite{CPRY2019}). Based on above the level-dependent scale functions it has been shown in \cite{CPRY2019},  for $\mathbb{W}^{(q)}(x)=\mathbb{W}^{(q)}(x,0)$, $0 \leq x \leq a$ and a Borel set $\mathfrak B \subseteq \mathbb R$, that
\begin{align}
	\begin{split}	\Ex_x\bigl[e^{-q \kappa_a^{+}} \mathbf{1}_{(\kappa_a^{+}<\kappa_0^{-})}\bigr]
&=\frac{\mathbb{W}^{(q)}(x)}{\mathbb{W}^{(q)}(a)}, \\
 	\int_0^{\infty} e^{-qt} \mathbb P_x(U_t \in \mathfrak{B},t<\kappa_a^{+}\wedge \kappa_0^{-})\dd dt
&=\int_{\mathfrak{B} \cap(0, a)} \Xi_\phi(z)^{-1}\Bigl(\frac{\mathbb{W}^{(q)}(x)}{\mathbb{W}^{(q)}(a)}\mathbb{W}^{(q)}(a; z)-\mathbb{W}^{(q)}(x; z)\Bigr)\mathrm{d}z.
\label{eqn:twoSidedALevel}
\end{split}
\end{align}
\section{Omega-killed level-dependent scale functions}
\label{killed scale}
In this section we introduce  the $\omega$-killed scale functions for the level-dependent process $U$  as part of solutions of renewal equations.
All fluctuation identities derived below are based on a new family of scale functions,
$\mathpzc{W}^{(\omega)}(x;y)$ and $\mathpzc{Z}^{(\omega)}(x;y)$, $x,\; y\in \R$,
which we shall call as \textit{$\omega$-killed level-dependent scale functions},
and are defined as the unique solutions to the following integral
equations
\begin{align}
	\mathpzc{W}^{(\omega)}(x, y) &= \mathbb{W}^{(q)}(x;y) + \int_{y}^{x}\mathbb{W}^{(q)}(x;z)\Xi_\phi(z)^{-1}(\omega(z)-q)\mathpzc{W}^{(\omega)}(z, y)\mathrm{d}z,
	\label{eqn:shiftedOmegaLevelScaleWq} \\
	\mathpzc{Z}^{(\omega)}(x,y) &=  \mathbb{Z}^{(q)}(x;y) + \int_{y}^{x}\mathbb{W}^{(q)}(x; z)\Xi_\phi(z)^{-1}(\omega(z)-q)\mathpzc{Z}^{(\omega)}(z, y)\mathrm{d}z,
	\label{eqn:shiftedOmegaLevelScaleZq}
\end{align}
where $\mathbb{W}^{(q)}(x;y)$, $\mathbb{Z}^{(q)}(x;y)$  and $\Xi_\phi(x)$ as in Eqs \eqref{eqn:shiftedLevelScaleW}-\eqref{eqn:Xi}
with $\mathpzc{W}^{(\omega)}(y, y)  = W^{(q)}(0)$
and $\mathpzc{Z}^{(\omega)}(y,y)  = 1$.
When $y=q=0$, we denote $	\mathpzc{W}^{(\omega)}(x, 0):=\mathpzc{W}^{(\omega)}(x)$ and $	\mathpzc{Z}^{(\omega)}(x,0):=	\mathpzc{Z}^{(\omega)}(x)$ and the above Eqs are reduced to
	\begin{align}
		\mathpzc{W}^{(\omega)}(x) &= \mathbb{W}(x) + \int_{0}^{x}\mathbb{W}(x;y)\Xi_\phi(y)^{-1}\omega(y)\mathpzc{W}^{(\omega)}(y)\mathrm{d}y,
		\label{eqn:omegaLevelScaleW} \\
		\mathpzc{Z}^{(\omega)}(x) &= 1 + \int_{0}^{x}\mathbb{W}(x; y)\Xi_\phi(y)^{-1}\omega(y)\mathpzc{Z}^{(\omega)}(y)\mathrm{d}y,
		\label{eqn:omegaLevelScaleZ}
	\end{align}
since $\mathbb{Z}(x;0)=Z(x)=1$.
In the remaining of the paper we shall use extensively integral equations of the form in Eqs \eqref{eqn:shiftedOmegaLevelScaleWq}-\eqref{eqn:omegaLevelScaleZ}.
The existence and uniqueness of a solution of such  integral equations it is given by the proposition below.
\begin{proposition}
		\label{lem:ExistenceOLS}
	Let $h(x, y)$ be a locally bounded function on $\R ^2$ and $\omega(\cdot)$ is a locally bounded function on
	$\R$. The following equation
		\begin{equation}
		H^{(\omega)}(x, y) = h(x, y) + \int_{y}^{x}\mathbb{W}^{(q)}(x;z)\Xi_\phi(z)^{-1}(\omega(z) - q)H^{(\omega)}(z, y)\mathrm{d}z,
		\label{eqn:EU1}
	\end{equation}
	admits a unique locally bounded solution on $\R^2$ and for  $x \leq y$  satisfies $H^{(\omega)}(x, y) = h(x, y)$.
\end{proposition}
\begin{proof}
	We start with existence of the solution. We note  that $\mathpzc{W}^{(q)}(x,y)$ is well-defined  and satisfies Eq. \eqref{eqn:shiftedOmegaLevelScaleWq}.
	For any $h(x,y)$, we define
		\begin{equation}
			H^{(\omega)}(x, y) := h(x, y) + \int_{y}^{x}\mathpzc{W}^{(\omega)}(z, y)\Xi_\phi(z)^{-1}(\omega(z) - q)h(z, y)\mathrm{d}z.
			\label{eqn:Hdef}
		\end{equation}
	Let $x > y$, based on the above equation, using Fubini's theorem and \eqref{eqn:shiftedOmegaLevelScaleWq}, we have
		\begin{align*}
	&\int_{y}^{x}\mathbb{W}^{(q)}(x;z)\Xi_\phi(z)^{-1}(\omega(z)-q)H^{(\omega)}(z, y)\mathrm{d}z \\
			&  =
			\int_{y}^{x}\mathbb{W}^{(q)}(x;z)\Xi_\phi(z)^{-1}(\omega(z)-q)h(z, y)\mathrm{d}z\\
			&\quad +\int_{y}^{x}\int_{u}^{x}\mathbb{W}^{(q)}(x;z)\Xi_\phi(z)^{-1}(\omega(z)-q)\mathpzc{W}^{(\omega)}(z, u)\Xi_\phi(u)^{-1}(\omega(u) - q)h(u, y)\mathrm{d}z\mathrm{d}u\\
			& =\int_{y}^{x}\mathpzc{W}^{(\omega)}(x, u)\Xi_\phi(u)^{-1}(\omega(u) - q)h(u, y)\mathrm{d}u\\
			& = H^{(\omega)}(x, y) - h(x, y).
					\end{align*}
To show the uniqueness of~\eqref{eqn:EU1}, suppose that there exist two solutions $H_1^{(\omega)}$ and $H_2^{(\omega)}$, then
	\( \widetilde{H}^{(\omega)} = H_1^{(\omega)} - H_2^{(\omega)}\) is solution to
	$\lvert \widetilde{H}^{(\omega)}(x, y)\rvert = \int_{y}^{x}\mathbb{W}^{(q)}(x;z)\Xi_\phi(z)^{-1} \lvert\omega(z) - q\rvert \lvert\widetilde{H}^{(\omega)}(z, y)\rvert \mathrm{d}z$.  Under the assumption that the kernel is integrable,  Gronwall's  lemma   yields that  $\widetilde{H}^{(\omega)}(x, y) \equiv 0$.
\end{proof}
\noindent The next proposition can be used to show the Volterra equations  \eqref{eqn:shiftedOmegaLevelScaleWq}-\eqref{eqn:shiftedOmegaLevelScaleZq} are true.
	Consider the general version of Eqs \eqref{eqn:omegaLevelScaleW}-\eqref{eqn:omegaLevelScaleZ},
\begin{align}\begin{split}\mathpzc{W}^{(\omega)}(x,y) &= \mathbb{W}(x;y) + \int_{y}^{x}\mathbb{W}(x;z)\Xi_\phi(z)^{-1}\omega(z)\mathpzc{W}^{(\omega)}(z,y)\mathrm{d}z,\label{Eq1}\\
		\mathpzc{Z}^{(\omega)}(x,y) &= 1 + \int_{y}^{x}\mathbb{W}(x; y)\Xi_\phi(z)^{-1}\omega(z)\mathpzc{Z}^{(\omega)}(z,y)\mathrm{d}z.
		\end{split}
\end{align}
 Then, the following proposition holds.
\begin{proposition} 	\label{pro2w} Let $(\mathpzc{W}^{(\omega_1)}, \mathpzc{Z}^{(\omega_1)})$ and $(\mathpzc{W}^{(\omega_2)}, \mathpzc{Z}^{(\omega_2)})$ be  of the form of Eq. \eqref{Eq1} with respect to $\omega_1(\cdot)\geq 0$ and $\omega_2(\cdot)\geq 0$, respectively. Then, for $x$, $y\in \mathbb R$,
	\begin{align*}\mathpzc{W}^{(\omega_1)}(x,y)-\mathpzc{W}^{(\omega_2)}(x,y) &=  \int_{y}^{x}\mathpzc{W}^{(\omega_1)}(x,z)\Xi_\phi(z)^{-1}\bigl(\omega_1(z)-\omega_2(z)\bigr)\mathpzc{W}^{(\omega_2)}(z,y)\mathrm{d}z,\\
\mathpzc{Z}^{(\omega_1)}(x,y)-\mathpzc{Z}^{(\omega_2)}(x,y) &=  \int_{y}^{x}\mathpzc{Z}^{(\omega_1)}(x,z)\Xi_\phi(z)^{-1}\bigl(\omega_1(z)-\omega_2(z)\bigr)\mathpzc{Z}^{(\omega_2)}(z,y)\mathrm{d}z. 	\end{align*}
	\end{proposition}
\begin{proof}
	Using Eq. \eqref{eqn:Hdef}, the Volterra eq. for $\mathpzc{W}^{(\omega)}(x,y)$ has solution of the form
\begin{equation*}
	\mathpzc{W}^{(\omega)}(x,y)=\mathbb{W}(x;y)+\int_y^x\mathpzc{W}^{(\omega)}(x,z)\Xi_\phi(z)^{-1}\omega(z)\mathbb{W}(z;y)\mathrm  dz.
	\end{equation*}	
Applying the above eq.  in the first and third   line below and Eq. \eqref{Eq1}   in second  line below, we have
\begin{align*}
&	\int_{\mathbb R}	\mathpzc{W}^{(\omega_1)}(x,z)\Xi_\phi(z)^{-1}\omega_2(z)\mathpzc{W}^{(\omega_2)}(z,y)\mathrm d z\\
	&=	\int_{\mathbb R}	\Bigl[\mathbb{W}(x;z) + \int_{\mathbb R}\mathpzc{W}^{(\omega_1)}(x,u)\Xi_\phi(u)^{-1}\omega_1(u)\mathbb{W}(u;z)\mathrm{d}u\Bigr]\Xi_\phi(z)^{-1}\omega_2(z)\mathpzc{W}^{(\omega_2)}(z,y)\mathrm d z\\
	&=\mathpzc{W}^{(\omega_2)}(x,y)-\mathbb{W}(x;y)+\int_{\mathbb R}\mathpzc{W}^{(\omega_1)}(x,u)\Xi_\phi(u)^{-1}\omega_1(u)\bigl[\mathpzc{W}^{(\omega_2)}(u,y)-\mathbb{W}(u;y)\bigr]\mathrm du
	\\
	&=\mathpzc{W}^{(\omega_2)}(x,y)-\mathpzc{W}^{(\omega_1)}(x,y)+\int_{\mathbb R}\mathpzc{W}^{(\omega_1)}(x,z)\Xi_\phi(z)^{-1}\omega_1(z)\mathpzc{W}^{(\omega_2)}(z,y)\dd dz,	\end{align*}
which proves the first eq. of the proposition.  The second eq.  can be proved similarly.
\end{proof}
\noindent   Using Proposition \ref{pro2w}  and repeating the reasoning of Proposition 2 in \cite{BoliXiaowen2018}, Eq.  \eqref{eqn:shiftedOmegaLevelScaleWq} follows immediately. Same line of logic can be used to show Eq. \eqref{eqn:shiftedOmegaLevelScaleZq}.
\subsection{Two-sided exit times  and potential measure}
In this section we shall use the $\omega$-killed level-dependent scale functions defined before to derive explicit expressions for the two-sided exit (above and below) and the killed potential measure.
\begin{theorem}	\label{thm:twoSidedAboveLevelw}
	For $0 \leq x \leq a$, we have
	\begin{equation}
		\label{eqn:twoSidedALevelw}
		\mathpzc{A}(x, a) = \Ex_x\Bigl[ e^{-\int_{0}^{\kappa_a^+}\omega(U_t)\mathrm{d}t}\mathbf{1}_{(\kappa_a^+< \kappa_0^-)}\Bigr]
		=\frac{\mathpzc{W}^{(\omega)}(x)}{\mathpzc{W}^{(\omega)}(a)}.
	\end{equation}
\end{theorem}

\begin{proof}
	 For
	$x \geq 0$, applying strong Markov property
	and using the fact that $U$ has no positive jumps gives that
		\begin{equation}
			\mathpzc{A}(x, a) = \mathpzc{A}(x, y)\mathpzc{A}(y, a), \quad \mbox{for all} \quad
			0 \leq x \leq y \leq a.
			\label{eqn:Arelation}
		\end{equation}
	Now, similarly to  \cite{boli2018}, we let $\omega$ be a bounded function and $\lambda$ be its arbitrary bound.  Define ${E} = \{{E}_t, \; t \geq 0\}$ to be a Poisson point process with measure $\frac{1}{\lambda}\mathbf{1}_{(0,\lambda]}(y)\mathrm{d}y \lambda \mathrm{d}t$. That is that  $E = \{(T_i, M_i), \; i = 1,2, \dots\}$ is a marked Poisson process with intensity $\lambda$ at jump epochs $T_i$ and (independent) marks $M_i$ being uniformly distributed on $[0,\lambda]$.
Therefore,
\[			\mathpzc{A}(x, a)
			= \Prob_x\bigl( \#_i\{ M_i < \omega(U_{T_i}) \; \mbox{for all} \; T_i \leq \kappa^+_a  \;
			\mbox{and} \; \{\kappa^+_a < \kappa^-_{0}\} \}=0 \bigr).
\]
Using Eq. \eqref{eqn:twoSidedALevel} and the Markov property of  $X$, we have
		\begin{align*}
			\mathpzc{A}(x, a)&=\mathbb P_x(T_1>\tau_a^+,\tau_a^+<\tau_0^-)+\Ex_x \bigl[ \mathpzc{A}(U_{T_1}, a)\mathbf{1}_{(T_1< \kappa^+_a  \wedge \kappa^-_{0}, \; M_1 >\omega(U_{T_1}))}\bigr]\\
			&= \frac{ \mathbb{W}^{(\lambda)}(x)}{ \mathbb{W}^{(\lambda)}(a)}
		+\int_{0}^{a} (\lambda -\omega(y))\mathpzc{A}(y, a)\Xi_\phi(y)^{-1}\Bigl(\frac{\mathbb{W}^{(\lambda)}(x)}{\mathbb{W}^{(\lambda)}(a)}\mathbb{W}^{(\lambda)}(a; y)-\mathbb{W}^{(\lambda)}(x; y)\Bigr)\mathrm{d}y.
		\end{align*}
Using Eq. \eqref{eqn:Arelation}, and re-arranging the above equation  gives
		\begin{align*}
		&	\mathpzc{A}(x, a) \Bigl(1 + \int_{0}^{x}\mathpzc{A}(y, x)(\lambda -\omega(y))\Xi_\phi(y)^{-1}\mathbb{W}^{(\lambda)}(x; y)\mathrm{d}y\Bigr)\\
		&	= \frac{ \mathbb{W}^{(\lambda)}(x)}{ \mathbb{W}^{(\lambda)}(a)}
		\Bigl( 1 + \int_{0}^{a}\mathpzc{A}(y, a)(\lambda -\omega(y))\Xi_\phi(y)^{-1}\mathbb{W}^{(\lambda)}(a; y)\mathrm{d}y\Bigr).
		\end{align*}
	Now, defining
		\begin{equation*}
			\mathpzc{W}^{(\omega)}(x) := \mathbb{W}^{(\lambda)}(x)\Bigl(1 +\int_{0}^{x}\mathpzc{A}(y, x)(\lambda -\omega(y))\Xi_\phi(y)^{-1}\mathbb{W}^{(\lambda)}(x; y)\mathrm{d}y\Bigr)^{-1},
			\label{eqn:OLdefW}
		\end{equation*}
it yields the required identity. Further, substituting  the form of  Eq. \eqref{eqn:twoSidedALevelw} into the above equation and noticing that $\mathpzc{A}(x, a)$ does not depend on $\lambda$, thus the above eq. should  be true for  $\lambda = 0$, completes the proof.
\end{proof}
\noindent In the next theorem we derive an expression for the $\omega$-type resolvent. To do  this we recall that 	for   a Markov process  $Z = \{Z_t, t\geq 0\} $ with killing time $\zeta$, and  $f$ a  non-negative bounded continuous function on $\R$  such that $\int_{0}^{\infty}\mathbb E_x[f(Z_t)1_{(t<\zeta)}] \dd dt < \infty$, the $\omega$-type resolvent, $K^{(\omega)}$, is given by
	\begin{equation*}
		K^{(\omega)}f(x) := \int_{0}^{\infty}Q_t^{(\omega)}f(x)\mathrm{d}t, \quad \text{with}\quad 		Q_t^{(\omega)}f(x) :=  \Ex_x\Bigl[
		{e^{-\int_{0}^{t}\omega(Z_s)\mathrm{d}s}} f(Z_t)\mathbf{1}_{(t < \zeta)}
		\Bigr].
	\end{equation*}
Then, (see Lemma 4.1 in \cite{boli2018})   $K^{(\omega)}f(x)$ is finite and satisfies the  equation
\begin{equation}
	K^{(\omega)}f(x) =
	K^{(0)}\bigl(f- \omega(\cdot) K^{(\omega)}f\bigr)(x).
	\label{eqn:resolventK}
\end{equation}
\begin{theorem}	\label{thm:wLevelResolventU}
	For $x, y \in [0, a]$, the $\omega$-type resolvent killed on exiting $[0,a]$ is given by
		\begin{align*}
						\label{eqn:wLevelResolvent}
		 \int_{0}^{\infty}\Ex_x\bigl[ e^{-\int_{0}^{t}\omega(U_s)\mathrm{d}s} \mathbf{1}_{(U_t \in \mathrm{d}y, \; t<\kappa_0^-\wedge\kappa_a^+ )}\bigr] \dd dt
			=\Xi_\phi(y)^{-1}\Bigl( \frac{\mathpzc{W}^{(\omega)}(x)}{\mathpzc{W}^{(\omega)}(a)}\mathpzc{W}^{(\omega)}(a,y)-
			\mathpzc{W}^{(\omega)}(x,y)\Bigr)\mathrm{d}y.
		\end{align*}
\end{theorem}
\begin{proof}
Following the line of logic in \cite{boli2018}, applying Eq. \eqref{eqn:resolventK} and  using Eq.  \eqref{eqn:twoSidedALevel}, it yields that 	
		\begin{align*}
			\mathpzc{U}^{(\omega)}f(x) &:= \int_{0}^{\infty}\Ex_x \Bigl[f(U_t)e^{-\int_{0}^{t}\omega(U_s)\mathrm{d}s}\mathbf{1}_{(t<\kappa_0^-\wedge\kappa_a^+ )}\Bigr]\mathrm{d}t\\
						&= \int_{0}^{a}\bigl(f(y)-\omega(y)\mathpzc{U}^{(\omega)}f(y)\bigr)\int_{0}^{\infty}\Prob_x(U_t \in \mathrm{d}y,\;t<\kappa_0^-\wedge\kappa_a^+ )\mathrm{d}t\\
			& = a_U\mathbb{W}(x) - \int_{0}^{x}\mathbb{W}(x; y)\Xi_\phi(y)^{-1}f(y)\mathrm{d}y
			+\int_{0}^{x}\mathbb{W}(x; y)\Xi_\phi(y)^{-1}\omega(y)\mathpzc{U}^{(\omega)}f(y)\mathrm{d}y.			
		\end{align*}
with
	$a_U := \int_{0}^{a}\bigl(f(y)-\omega(y)\mathpzc{U}^{(\omega)}f(y)\bigr)\Xi_\phi(y)^{-1}
		\frac{\mathbb{W}(a; y)}{\mathbb{W}(a)}\mathrm{d}y$.
		Now, we define the operator
		\begin{equation}
			\mathpzc{R}^{(\omega)}f(x) := \int_{0}^{x}f(y)\Xi_\phi(y)^{-1}\mathpzc{W}^{(\omega)}(x,y)\mathrm{d}y, \quad
			x > 0,
			\label{eqn:defR}
		\end{equation}
	with $\mathpzc{R}^{(\omega)}f(x) = 0$ for $x \leq 0$. Using Eq. \eqref{eqn:shiftedOmegaLevelScaleWq},
	we obtain that
		\begin{align*}
			\mathpzc{R}^{(\omega)}f(x) &=\int_{0}^{x}\mathbb{W}(x;y)\Xi_\phi(y)^{-1}f(y)\mathrm{d}y +\int_{0}^{x}\mathbb{W}(x; y)\Xi_\phi(y)^{-1}\omega(y)\mathpzc{R}^{(\omega)}f(y)\mathrm{d}y.
			\label{eqn:operatorR}
		\end{align*}
	Since the form of the above Voterra-type eqs are uniquely defined,
	we conclude that
		\begin{equation*}
			\mathpzc{U}^{(\omega)}f(x) = a_U\mathpzc{W}^{(\omega)}(x) - \mathpzc{R}^{(\omega)}f(x),
			\label{eqn:resolventRelation}
		\end{equation*}
which along with the  boundary condition $\mathpzc{U}^{(\omega)}f(a) = 0$, gives that
		\begin{equation}
			\mathpzc{U}^{(\omega)}f(x) = \int_{0}^{\infty}f(y)\Xi_\phi(y)^{-1}
			\Bigl(\frac{\mathpzc{W}^{(\omega)}(x)}{\mathpzc{W}^{(\omega)}(a)} \mathpzc{W}^{(\omega)}(a,y) - \mathpzc{W}^{(\omega)}(x,y)\Bigr)\mathrm{d}y,
			\label{eqn:Ufx}
		\end{equation}
and thus the proof is completed.
\end{proof}
\noindent Now, using Theorem \ref{thm:wLevelResolventU}, we get the following identity for the two-sided exit from below problem.
\begin{theorem}	\label{thm:twoSidedBelowLevelw}
	For $x \leq a$, we have
	\begin{equation*}
		\mathpzc{B}(x, a) = \Ex_x\Bigl[ e^{-\int_{0}^{\kappa_0^-}\omega(U_t)\mathrm{d}t}\mathbf{1}_{( \kappa_0^- <\kappa_a^+)}\Bigr]
		=\mathpzc{Z}^{(\omega)}(x) -\frac{\mathpzc{W}^{(\omega)}(x)}{\mathpzc{W}^{(\omega)}(a)} \mathpzc{Z}^{(\omega)}(a).
		\label{eqn:twoSidedBLevelw}
	\end{equation*}
\end{theorem}
\begin{proof}
	First note that recalling the definition of  $\mathpzc{Z}^{(\omega)}(x)$ and  using Proposition \ref{lem:ExistenceOLS}, we have that
	\begin{equation}
		\label{solZx}
	\mathpzc{Z}^{(\omega)}(x)=1+\int_0^x	\mathpzc W ^{(\omega)}(x,z)	\Xi_\phi(z)^{-1}\omega(z)\dd dz.
		\end{equation}
	Further, using Eq. \eqref{eqn:Ufx}, it yields that
	\begin{align*}
			\mathpzc{U}^{(\omega)}\omega(x) &= \int_{0}^{a}\omega(y)\Xi_\phi(y)^{-1}
		\Bigl(\frac{\mathpzc{W}^{(\omega)}(x)}{\mathpzc{W}^{(\omega)}(a)} \mathpzc{W}^{(\omega)}(a,y) - \mathpzc{W}^{(\omega)}(x,y)\Bigr)\mathrm{d}y\\
		&=\frac{\mathpzc{W}^{(\omega)}(x)}{\mathpzc{W}^{(\omega)}(a)}\bigl(\mathpzc{Z}^{(\omega)}(a)-1\bigr)-\bigl(\mathpzc{Z}^{(\omega)}(x)-1\bigr),
	\end{align*}
where the last eq. follows from Eq. \eqref{solZx}. On the other hand, using Theorem \ref{thm:twoSidedAboveLevelw}, noticing that
\begin{equation*}
		\mathpzc{U}^{(\omega)}\omega(x) =1-\mathbb E_x\Bigl[e^{-\int_0^{ \kappa_0^-\wedge\kappa_a^+  }\omega(U_t)\dd dt }\Bigr]=1-\frac{\mathpzc{W}^{(\omega)}(x)}{\mathpzc{W}^{(\omega)}(a)}-	\mathpzc{B}(x, a) ,
	\end{equation*}
and combining the above two eqs.,  the result follows.
	\end{proof}
\begin{remark}\label{rem1}
 We note that the exit identities in Theorems \ref{thm:twoSidedAboveLevelw} and \ref{thm:twoSidedBelowLevelw} can be generalised  for any interval $[y,z]$, i.e. for  $y \leq x \leq z$, it holds that
			\begin{align*}
				\Ex_x\Bigl[ e^{-\int_{0}^{\kappa_z^+}\omega(U_t)\mathrm{d}t}\mathbf{1}_{(\kappa_z^+< \kappa_y^-)}\Bigr]
				&=\frac{\mathpzc{W}^{(\omega)}(x,y)}{\mathpzc{W}^{(\omega)}(z,y)},
		\\
				\Ex_x\Bigl[ e^{-\int_{0}^{\kappa_y^-}\omega(U_t)\mathrm{d}t}\mathbf{1}_{( \kappa_y^- <\kappa_z^+)}\Bigr]
				&=\mathpzc{Z}^{(\omega)}(x,y) -\frac{\mathpzc{W}^{(\omega)}(x,y)}{\mathpzc{W}^{(\omega)}(z,y)} \mathpzc{Z}^{(\omega)}(z,y).
			\end{align*}
	\end{remark}
\subsection{One-sided upward exit problem}
\label{subsec:OneSided}
 To solve one-sided upwards exit problem,
we assume that
	$	\omega(x) = p$, for all $x \leq 0$,
and we define the function $\mathpzc{H}^{(\omega)}(x)$ on $\R$ satisfying the following integral equation
	\begin{equation}		\label{eqn:HRenewal}
		\mathpzc{H}^{(\omega)}(x) = u^{(p)}(x)
		+ \int_{0}^{x}\mathbb{W}^{(p)}(x;y)\Xi_\phi(y)^{-1}(\omega(y)-p)
		\mathpzc{H}^{(\omega)}(y)\mathrm{d}y,
	\end{equation}
where $u^{(p)}(x):= \lim_{\gamma \to \infty}\mathbb{W}^{(p)}(x; -\gamma)/W^{(p)}(\gamma)$ satisfies
\[u^{(p)}(x)=e^{\phi_px}+\int_{-\gamma}^x  W^{(q)}(x-y) \phi(y)u^{({p}) \prime}(y) \mathrm{d}y,\]
with $\phi_p$ to be right inverse $\phi_p=\sup\{\theta\geq 0:\psi(\theta)=p\}$ and ${u^{(p)\prime}}$ to be  the unique solution to
	\begin{equation*}
		u^{(p)\prime}(x)=\Xi_\phi(x)^{-1} \phi_p e^{\phi_p x}+\int_{-\gamma}^x \Xi_\phi(x)^{-1} \phi(y) W^{(p)\prime}(x-y) u^{(p)\prime}(y)\mathrm{d} y.
	\end{equation*}

\begin{theorem}	\label{thm:OneSidedUpward}
Assume that $\omega(x) = p$, for all $x \leq 0$. Then,  for $x \leq a$, we have that \upshape
		\begin{equation*}
			\Ex_x\Bigl[e^{-\int_0^{\kappa_a^{+}} \omega(U_s)\mathrm{d}s}\mathbf{1}_{(\kappa_a^{+}<\infty)}\Bigr]
			= \frac{\mathpzc{H}^{(\omega)}(x)}{\mathpzc{H}^{(\omega)}(a)},
			\label{eqn:OneSidedUpward}
		\end{equation*}
	with the corresponding resolvent for $x, y \leq a$, to be given by
		\begin{equation*}
		 \int_0^{\infty}\Ex_x\Bigl[e^{-\int_0^t \omega(U_s)\mathrm{d}s} \mathbf{1}_{(U_t \in \mathrm{d}y, \; t<\kappa_a^{+})}\Bigr] \mathrm{d}t
			=\Xi_\phi(y)^{-1}\Bigl(\frac{\mathpzc{H}^{(\omega)}(x)}{\mathpzc{H}^{(\omega)}(a)} \mathpzc{W}^{(\omega)}(a, y)-\mathpzc{W}^{(\omega)}(x, y)\Bigr) \mathrm{d}y,
			\label{eqn:ResolventUpward}
		\end{equation*}
where $\mathpzc{H}^{(\omega)}$ is given in~\eqref{eqn:HRenewal}.
\end{theorem}
\begin{proof}
For the first eq. of the theorem, it suffices to prove
	\begin{equation*}
	\lim_{\gamma \to \infty}\frac{\mathpzc{W}^{(\omega)}(x,-\gamma)}{\mathpzc{W}^{(\omega)}(a, -\gamma)}
	= \frac{\mathpzc{H}^{(\omega)}(x)}{\mathpzc{H}^{(\omega)}(a)}.
	\label{eqn:Homega}
\end{equation*}
To do this,  note that using   Remark \ref{rem1} it holds that
	\begin{equation}
		\mathpzc{W}^{(\omega)}(x, -\gamma) = \mathpzc{W}^{(\omega)}(0,-\gamma)	\Bigl(\Ex\Bigl[ e^{-\int_{0}^{\kappa_x^+}\omega(U_t)\mathrm{d}t}\mathbf{1}_{(\kappa_x^+< \kappa_{-\gamma}^-)}\Bigr] \Bigr)^{-1},
		\label{eqn:expectH}
	\end{equation}
and thus the above limit should be identified.
We note  first that
		\begin{align*}
			\mathpzc{W}^{(\omega)}(x, -\gamma) &= \mathbb{W}^{(p)}(x; -\gamma) + \int_{0}^{x}\mathbb{W}^{(p)}(x;z)\Xi_\phi(z)^{-1}(\omega(z)-p)\mathpzc{W}^{(\omega)}(z, -\gamma)\mathrm{d}z,
			\label{eqn:gammaW}
		\end{align*}
and thus for $x \in [-\gamma, 0]$,
	 we have that
		\[\lim_{\gamma \to \infty}\frac{\mathpzc{W}^{(\omega)}(x, -\gamma)}{W^{(p)}(\gamma)} =
		\lim_{\gamma \to \infty}\frac{\mathbb{W}^{(p)}(x; -\gamma)}{W^{(p)}(\gamma)} =u^{(p)}(x),\]
where $u^{(p)}(x)$ is  the same form of Lemma 28 in \cite{CPRY2019}, satisfying the above integral eqs.
Clearly, for $x = 0$, it yields that
$\lim_{\gamma \to \infty}\mathpzc{W}^{(\omega)}(0, -\gamma)/W^{(p)}(\gamma)= u^{(p)}(0) = 1.$

\noindent 	Moreover, since the expectation in Eq. \eqref{eqn:expectH} is (monotone) increasing with respect to $\gamma $, then the following limit is well defined and finite for $x \geq - \gamma$, which gives
		\begin{equation*}
			\lim_{\gamma \to \infty}\frac{\mathpzc{W}^{(\omega)}(x, -\gamma)}{W^{(p)}(\gamma)} = \Bigl(\Ex\Bigl[ e^{-\int_{0}^{\kappa_x^+}\omega(U_t)\mathrm{d}t}\mathbf{1}_{(\kappa_x^+< \infty)}\Bigr] \Bigr)^{-1}.
		\end{equation*}
Defining
$			\mathpzc{H}^{(\omega)}(x)
			:= \lim_{\gamma \to \infty}\mathpzc{W}^{(\omega)}(x, -\gamma)/W^{(p)}(\gamma),
$
	completes the first part of the theorem.
		Furthermore, we show that
	the above-defined $\mathpzc{H}^{(\omega)}(x)$ satisfies Eq. \eqref{eqn:HRenewal}.
	Multiplying the above renewal eq.  with $1/ W^{(p)}(\gamma)$,  taking limit as $\gamma \to \infty$ and  using  dominated convergence theorem, we obtain the required equation.
	
	Finally, to prove the resolvent identity of the theorem, from  Theorem \ref{thm:wLevelResolventU}, we have  that
		\begin{align*}
			\int_{0}^{\infty}\Ex_x\bigl[ e^{-\int_{0}^{t}\omega(U_s)\mathrm{d}s} \mathbf{1}_{(U_t \in \mathrm{d}y, \; t<\kappa_{-\gamma}^-\wedge\kappa_a^+ )}\bigr]
			=\Xi_\phi(y)^{-1}\Bigl( \frac{\mathpzc{W}^{(\omega)}(x, -\gamma)}{\mathpzc{W}^{(\omega)}(a, -\gamma)}\mathpzc{W}^{(\omega)}(a,y)-
			\mathpzc{W}^{(\omega)}(x,y)\Bigr)\mathrm{d}y,
		\end{align*}
	from which by  taking limit of the both sides as $\gamma \to \infty$  and
	using Eq. \eqref{eqn:Homega} we recover the  required equation.
\end{proof}
\subsection{Omega  level-dependent risk model}
In this section, we demonstrate one application of previously obtained results to the omega model by computing the bankruptcy probability for
an omega-killed level-dependent L\'evy process $U$ with rate function $\phi$, under the assumption that $\Ex(X_1) = \psi'(0+)> 0$.
In  the omega model the insurance company
is allowed to continue its business, even with a negative surplus,  till bankruptcy,   as it was proposed in \cite{GSH2012} and  studied for the  L\'evy risk process  in \cite{AGS2011, GSH2012, LRZ2011, LRZ2014}. Extending  these  results, within the level-dependent Omega risk model,
bankruptcy happens either when the surplus process in $[-d, 0]$ and the bankruptcy rate is a function of the current level, or
occur if the surplus process falls below some fixed level  $-d < 0$.

We consider the omega-killed level-dependent L\'evy process as the surplus process of an insurer and $\omega_d(x)=
	\omega(x) \mathbf 1_{(x\in [-d, 0])}$.
Clearly,  $x < 0$, the quantity of $\omega(x)dt$ describes the
	probability of bankruptcy within an infinitesimal time units, $dt$ and
 the bankruptcy probability is defined by
		\begin{equation*}
			\mathcal{P}(x) = 1 - \Ex_x \bigl[e^{-\int_0^{\infty} \omega\left(U_s\right)\mathrm{d}s} \mathbf{1}_{(\tau_{-d}^{-}=\infty)}\bigr].
		\end{equation*}
	\begin{proposition}
		\label{prop:bankrupcty}
		Assume that $x \in [-d, \infty)$.  If $\Ex(X_1) > 0$, and $\int_{-d}^{\infty} \phi(x) \mathbb{W^{\prime}}(x;-d)\mathrm{d}x$ exists, then the bankruptcy probability is  given by
				\begin{equation*}
					\mathcal{P}(x) = 1 - a_{\mathpzc{W}^{-1}(\infty ,-d)}\mathpzc{W}^{(\omega)}(x, -d),
				\end{equation*}
			where
			\begin{linenomath*}
				\begin{equation*}
					a_{\mathpzc{W}^{-1}(\infty ,-d)^{-1}} = \frac{1}{A(-d) + \int_{-d}^{\infty}A(z)\Xi_\phi(z)^{-1}\omega(z)\mathpzc{W}^{(\omega)}(z)\mathrm{d}z },
				\end{equation*}
			\end{linenomath*}
			with \( A(x) := \frac{1}{\Ex(X_1)} \bigl( 1 + \int_{x}^{\infty} \phi(y) \mathbb{W^{\prime}}(y;-d)\mathrm{d}y\bigr).\)
				Furthermore, $\mathcal{P}(x)$ satisfies the following integral equation
			\begin{linenomath}
				\begin{equation*}
					1-\mathcal{P}(x) = a_{\mathpzc{W}^{-1}(\infty ,-d)}\mathbb{W}(x; -d) +
					\int_{-d}^{x}\mathbb{W}(x; z)\Xi_\phi(z)^{-1}\omega_d(z)(1-\mathcal{P}(z))\mathrm{d}z.
				\end{equation*}
			\end{linenomath}
			Moreover, when $\int_{x}^{\infty} \phi(y) \mathbb{W^{\prime}}(y;-d)\mathrm{d}y = \infty$, it follows that $A(x) = \infty$, hence $\mathcal{P}(x) = 1$.
	\end{proposition}
	
	\begin{proof}
	Recalling  Remark \ref{rem1},
		for $x \in [-d, 0]$, we have
			\begin{equation*}
				\mathpzc{W}^{(\omega)}(x, -d) = \mathbb{W}(x; -d) + \int_{-d}^{x}\mathbb{W}(x; z)\Xi_\phi(z)^{-1}\omega(z)\mathpzc{W}^{(\omega)}(z, -d)\mathrm{d}z.	
			\end{equation*}
		Taking limit $a \to \infty$  and using the fact that $	\lim _{a \to \infty}\mathbb{W}(a; -d) = \bigl( 1 + \int_{-d}^{\infty}\phi(y)\mathbb{W^{\prime}}(y; -d)\mathrm{d}y\bigr)/\Ex(X_1)$ [see Proposition 33 in \cite{CPRY2019}], gives
			\begin{equation*}
				\lim _{a \to \infty}\mathpzc{W}^{(\omega)}(a; -d) =A(-d) +
				\int_{-d}^{\infty}A(z)\Xi_\phi(z)^{-1}\omega(z)\mathpzc{W}^{(\omega)}(z, -d)\mathrm{d}z,
			\end{equation*}
	from which we get the expression  for $a_{\mathpzc{W}^{-1}(\infty ,-d)}=[\lim _{a \to \infty}\mathpzc{W}^{(\omega)}(a; -d)]^{-1}$.	
	Similarly, one can show the result, for 	 $x > 0 $, after noticing that
	in the positive half-line $\omega_d(x) = 0$, and thus the  limits in the in the integral of $a_{\mathpzc{W}^{-1}(\infty ,-d)}$ change from $-d$ to $0$ (instead of $[-d,x]$). Finally, by noticing that
			\begin{align*}
				\mathcal{P}(x)
				&=  1 - a_{\mathpzc{W}(\infty ,-d)^{-1}}\mathbb{W}(x; -d) - \int_{-d}^{x}\mathbb{W}(x; z)\Xi_\phi(z)^{-1}\omega_b(z)a_{\mathpzc{W}(\infty ,-d)^{-1}}\mathpzc{W}^{(\omega)}(z, -d)\mathrm{d}z\\
				&= 1 - a_{\mathpzc{W}(\infty ,-d)^{-1}}\mathbb{W}(x; -d) -
				\int_{-d}^{x}\mathbb{W}(x; z)\Xi_\phi(z)^{-1}\omega_b(z)(1-\mathcal{P}(z))\mathrm{d}z,
			\end{align*}
	we get the integral eq. of the proposition.
	\end{proof}

\section{Reflected omega-killed level-dependent L\'evy process}
\label{reflected}
In this section we introduce and study  the  case of a reflected  level-dependent L\'evy  process. Potential measures and identities for the exit problem are also derived.   To model our problem, we assume that level-dependent process $U$ in Eq. \eqref{eqn:levelDependent}  is reflected at a lower barrier and thus it is well known that can be expressed as the difference between the underlying level-dependent process and its running infimum process, see \cite{B2009, boli2018, PY2017, PY2018, P2004}, among others. It is worth mentioning that such models  have applications to risk theory where the reflection (at  a lower boundary  0) models the capital injections required so as the surplus process of insurance firm stays positive, see \cite{AE2014, PY2018}.
\subsection{Reflected level-dependent L\'evy process}
Let $V=(V_t)_{t\geq 0}$ be the  reflected version of $U$, at a lower boundary 0,  given by
	\begin{equation}
		V_t = X_t + R_t - \int_{0}^{t}\phi(V_s)\mathrm{d}s,
		\label{eqn:SDELevelRef}
	\end{equation}
where $R_t$ is a non-decreasing and right continuous process that represents the cumulative amounts of modification up to $t$ that pushes upward the process in Eq. \eqref{eqn:levelDependent} when it attempts to go below 0.   For the existence of the solution of the above eq., we use the theory of the multi-refracted-reflected L\'evy processes given in Appendix \ref{multi-reftract}. To end this, we approximate  a general function
$\phi$ by a sequence of a rate function $(\phi_n)_{n\geq 1}$ and  thereby we are able to define a sequence
of reflected multi-refracted L\'evy processes, which  allows us to show the existence of the solution of \eqref{eqn:SDELevelRef} by a limiting argument. We point out that under the assumptions  for $(\phi_n)_{n\geq 1}$ in the proof of  Proposition \ref{prop:convergenceVn} below,  \(V_t\)  in Eq.  \eqref{eqn:SDELevelRef}, is the result of   a uniform convergence of the  the multi-refracted-reflected L\'evy processes given in Appendix \ref{multi-reftract} (see proof of Proposition \ref{prop:convergenceVn}).
\begin{proposition}	\label{prop:convergenceVn}
	Suppose that the locally Lipschitz continous rate function $\phi$ satisfies the assumption \textbf{[A]}. Then, there exists a solution to $V_t$ to the SDE given in Eq. ~\eqref{eqn:SDELevelRef}.
\end{proposition}
\begin{proof}
	We choose as in  \cite{CPRY2019},    $(\phi_{n})_{n\geq 1}$ to be a  non-decreasing approximating sequence for $\phi$, that
	satisfies the following
	conditions
	\begin{itemize}
		\item[](i).  \(\lim_{n \to \infty}\phi_n = \phi\) uniformly on the compact sets,
		\item[](ii). \(\phi_1(x) \leq \phi_2(x) \leq \cdots \leq \phi(x)\) for all $x \in \R$,
		\item[] (iii).  For each $n\geq 1$, the rate function \( \phi_n(x) = \sum_{j=1}^{m_n}\delta_j^n\mathbf{1}_{(x > b_j^n)}\) for some $m_n \in \N$, $0 < b_1^n < \cdots < b_{m_n}^n$ and $\delta_j^n > 0$ for
		each $j = 1, \ldots, m_n$.
	\end{itemize}
We note that $\phi_n$ is an $m_n$ multi-refracted rate function (which is different of the rate function in Section \ref{multi-reftract}). For the construction of such function $\phi_n$ that satisfy the above conditions, we refer to Remark 18 of \cite{CPRY2019}.

For each
$n\geq 1$, we consider $V_n$ as a solution to the following equation (see Proposition \ref{prop:ConvergenceVk})
	\[ V_n(t) = X_t + R_n(t)- \int_{0}^{t} \phi_{n}(V_n(s))\mathrm{d}s, \quad t\geq 0.\]
For $\phi_1(x) \leq \phi_2(x) \leq \ldots \leq \phi(x)$ for all
$x \in \R$, we have from  Lemma \ref{lem:monotonocityVn} that $V_1(t) \geq V_2(t) \geq \ldots$.
Now, fix an arbitrary $T > 0$. It follows after defining $\overline{X}_T := \sup_{0 \leq s \leq T}X_s$ and $\underline{X}_T := \inf_{0 \leq s \leq T}X_s$ for any $T>0$ that	
\begin{align*}
	V_n(t)
	&\leq \bigl|\overline{X}_T\bigl|+ R_T, \qquad 0 \leq t \leq T,
\end{align*}
as $R_t$ is right-continuous and non-decreasing.
\noindent On the other hand, since $\phi(x) \geq 0$ for all $x\in \R$,
	\begin{align*}
		V_n(t) &\geq  X_t -R_T - \phi\bigl(\bigl|\overline{X}_T\bigl|+R_T\bigr)T,   \quad 0 \leq t \leq T .
	\end{align*}
Hence, from the last two eqs., for $n \geq 1$ and  $0 \leq t \leq T$, we get
	\[\bigl|V_n(t) \bigr| \leq \bigl(\bigl| \overline{X}_T\bigr| \vee \bigl| \underline{X}_T\bigr| \bigr)+ R_T +\phi\bigl(\bigl|\overline{X}_T\bigl|+ R_T	\bigr)T := K_T. \]
Since by~Lemma \ref{lem:monotonocityVn}, the sequence $n \mapsto V_n$ is non-increasing  and
bounded below, we can define \(V(t) := \lim_{n \to \infty}V_n(t) \).
Further, by the uniform convergence of $\phi_{n} \to \phi$ on compact sets
and pointwise convergence of $V_n \to V$, we have
$\phi_n(V_n) \to \phi(V)$ pointwise.

\noindent  To show uniform convergence of $V_n$ to $V$ on compact time intervals,
we consider
	\begin{align*}
		V_{n}(t) - V_t
		&= R_n(t)-R_t
		 + \int_{0}^{t}\bigl(\phi(V_s)
		-\phi_n(V_n(s))\bigr)\mathrm{d}s\leq \int_{0}^{t}\bigl(\phi(V_s)
		-\phi_n(V_n(s)) \bigr)\mathrm{d}s,
	\end{align*}
which the last inequality follows from  Lemma \ref{lem:monotonocityVn}.
Then, using triangular inequality,  the assumption of $\phi$ is a locally Lipschitz function and Mean Value Theorem,  we have
	\begin{align*}
		|V_{n}(t) - V_t| &\leq \int_{0}^{t}\bigl|\phi(V_s)- \phi(V_n(s))\bigr|\mathrm{d}s +
		\int_{0}^{t}\bigl|\phi(V_n(s)) -\phi_n(V_n(s)) \bigr|\mathrm{d}s\\
		&\leq L_I\int_{0}^{t}\bigl|V_s- V_n(s)\bigr|\mathrm{d}s +
		T\sup_{s \in I} \bigl|\phi(s) -\phi_n(s)\bigr|,
	\end{align*}
where $I := [-K_T, K_T]$ and $L_I$  is the associated Lipschitz constant for $\phi$ over the interval $I$.
Thus, using Gronwall's Inequality  gives that
	\[\sup_{0 \leq t \leq T}|V_{n}(t) - V_t| \leq T\sup_{s \in I} \bigl|\phi(s) -\phi_n(s) \bigr|e^{L_IT}.\]
Since $\lim_{n\to \infty}\phi_{n} = \phi$ uniformly a.s., we conclude that $V_{n}(t)$
converges uniformly to $V(t)$ a.s. on compact time intervals.
\end{proof}
\subsection{Exit time and potential measure}
\label{killed reflection}
In this section we derive identities for the exit time upwards and the potential measure of  the $\omega$-killed  reflected level-dependent L\'evy process in Eq. \eqref{eqn:SDELevelRef}. Let  $	K_a$ be
the first passage time for the reflected level-dependent L\'evy process exits above a level $a \in \R$, given by
	\begin{equation*}
		K_a = \inf\{ t> 0: V_t > a\}, \; \mbox{ for } \; a > 0.
	\end{equation*}
\begin{theorem}	\label{thm:reflectedI_OmegaLevel}
	For $0 \leq x \leq a$, we have
	\begin{linenomath*}
		\begin{equation*}
			\mathpzc{C}(x, a) := \Ex_x\bigl[ e^{-\int_{0}^{K_a}\omega(V_t)\mathrm{d}t}\mathbf{1}_{(K_a < \infty)}\bigr]
			= \frac{\mathpzc{Z}^{(\omega)}(x)}{\mathpzc{Z}^{(\omega)}(a)}.
		\end{equation*}
	\end{linenomath*}
\end{theorem}

\begin{proof}
	For all $0\leq x \leq y\leq z$,  applying strong Markov property of $V_t$ and
	using the
	fact that $V_t$ has no positive jumps, we have
		\begin{equation}\label{rmp}
			\mathpzc{C}(x, z) = \mathpzc{C}(x, y)\mathpzc{C}(y, z), \quad \mbox{for all} \quad
			0\leq x\leq y\leq z.
		\end{equation}
	Moreover, for $x\leq a$, we have
	\begin{linenomath*}
		\begin{align*}
			1- \mathpzc{C}(x, a) &= \Ex_x\Bigl[ \int_{0}^{K_a}e^{-\int_{0}^{t}\omega(V_s)\mathrm{d}s}\omega(V_t)\mathrm{d}t\Bigr]= \int_{0}^{\infty}\Ex_x\Bigl[\omega(V_t)e^{-\int_{0}^{t}\omega(V_s)\mathrm{d}s}\mathbf{1}_{(t <K_a)}\Bigr]\mathrm{d}t,
		\end{align*}
	\end{linenomath*}
which after applying Eq. \eqref{eqn:resolventK} becomes
		\begin{align*}
			\mathpzc{C}(x, a) &=1- K^{(0)}\bigl(\omega(\cdot)-\omega(\cdot)K^{(\omega)}\omega(\cdot)\bigr)(x)\\
		&=	1- \int_{0}^{a}\omega(y)\mathpzc{C}(y,a)
			\int_{0}^{\infty}\Prob_x(V_t \in \mathrm{d}y, \; t < K_a)\mathrm{d}t.
		\end{align*}
Now, note that  $\int_{0}^{\infty}\Prob_x(V_t \in \mathrm{d}y, \; t < K_a)\mathrm d t$ is the resolvent of
	level-dependent L\'evy process reflected at its infimum when q= 0, which is derived in Appendix \ref{A2} and is given by  \eqref{eqn:resolventLevelRef}. Using this and Eq. \eqref{rmp}, we get that
		\begin{align*}
			\mathpzc{C}(x, a) &=
			1- \int_{0}^{a}\omega(y)\mathpzc{C}(y,a)\Xi_\phi(y)^{-1}
			\bigl(\mathbb{W}(a;y) - \mathbb{W}(x;y) \bigr)\mathrm{d}y \\
			&=  1- \int_{0}^{a}\mathbb{W}(a;y)\Xi_\phi(y)^{-1}\omega(y)\mathpzc{C}(y,a)\mathrm{d}y -
			\mathpzc{C}(x,a)\int_{0}^{x}\mathbb{W}(x;y)\Xi_\phi(y)^{-1}\omega(y)
			\mathpzc{C}(y,x)\mathrm{d}y,
		\end{align*}
or equivalently
		\begin{equation*}
			\mathpzc{C}(x, a)\bigl[ 1- \int_{0}^{x}\mathbb{W}(x;y)\Xi_\phi(y)^{-1}\omega(y)\mathpzc{C}(y,x)\mathrm{d}y \bigr]	
			= 1- \int_{0}^{a}\mathbb{W}(a;y)\Xi_\phi(y)^{-1}\omega(y)\mathpzc{C}(y,a)\mathrm{d}y.
		\end{equation*}
	Defining
		\begin{equation*}
			\mathpzc{Z}^{(\omega)}(x):= \Bigl(1- \int_{0}^{x}\mathbb{W}(x;y)\Xi_\phi(y)^{-1}\omega(y)\mathpzc{C}(y,x)\mathrm{d}y\Bigr)^{-1},
			\label{eqn:defZomegalevel}
		\end{equation*}
	gives the  required identity.
Finally, using the fact that  $\mathpzc{C}(y, x) =\mathpzc{Z}^{(\omega)}(y)/\mathpzc{Z}^{(\omega)}(x)$ in the above eq., produces~\eqref{eqn:omegaLevelScaleZ}.
\end{proof}
\noindent In the next theorem we derive an identity for the $\omega$-killed  resolvent measure of $V$.
\begin{theorem}
	For $x,y \in [0, a]$, we have
		\begin{align*}
			\mathpzc{L}^{(\omega)}(x, \mathrm{d}y) :=\int_0^{\infty} \Ex_x\Bigl[e^{-\int_0^t \omega(V_s)\mathrm{d}s}\mathbf{1}_{(V_t \in \mathrm{d}y, \; t<K_a)}\Bigr]\mathrm{d}t,
		\end{align*}
	which is a measure absolutely continuous with respect to the Lebesgue measure and has density
		\[\mathpzc{l}^{(\omega)}(x, y)=\Xi_\phi(y)^{-1} \Bigl(\frac{\mathpzc{Z}^{(\omega)}(x)}{\mathpzc{Z}^{(\omega)}(a)} \mathpzc{W}^{(\omega)}(a, y)-\mathpzc{W}^{(\omega)}(x, y)\Bigr).\]
\end{theorem}
\begin{proof}
	In this proof uses the same line of logic as Theorem \ref{thm:wLevelResolventU}. Let $f$ is non-negative bounded continuous function. Applying Eq. \eqref{eqn:resolventK}, we have
		\begin{align}
			\mathpzc{L}^{(\omega)}f(x) &:= \int_{0}^{\infty}\Ex_x \Bigl[f(V_t)e^{-\int_{0}^{t}\omega(V_s)\mathrm{d}s}\mathbf{1}_{(t<K_a)}\Bigr]\mathrm{d}t\nonumber\\
			&=\mathpzc{L}^{(0)}\bigl(f-\omega (\cdot)\mathpzc{L}^{(\omega)f}\bigr)(x)\nonumber\\
			&=\int_{0}^{\infty}\int_{0}^{a}\bigl(f(y)-\omega(y)\mathpzc{L}^{(\omega)}f(y)\bigr)\Prob_x(V_t \in \mathrm{d}y,\;t<K_a )\mathrm{d}t\nonumber\\
&=\int_{0}^{a}\mathbb{W}(a; y)\Xi_\phi(y)^{-1}\bigl(f(y)-\omega(y)\mathpzc{L}^{(\omega)}f(y)\bigr) \mathrm{d}y
- \int_{0}^{a}\mathbb{W}(x; y)\Xi_\phi(y)^{-1}f(y)\mathrm{d}y\nonumber\\
&\quad + \int_{0}^{a}\mathbb{W}(x; y)\Xi_\phi(y)^{-1} \omega(y)\mathpzc{L}^{(\omega)}f(y)\mathrm{d}y.
\label{eqn:L}
		\end{align}
	where the last equality follows from~Theorem \ref{thm:resolventLevelRef} in Appendix \ref{A2}.
	
	\noindent Defining
		\begin{equation*}
			a_L:= \int_{0}^{a}\mathbb{W}(a; y)\Xi_\phi(y)^{-1}\bigl(f(y)-\omega(y)\mathpzc{L}^{(\omega)}f(y)\bigr) \mathrm{d}y,
		\end{equation*}
	and recalling $\mathpzc{R}^{(\omega)}(x)$ given in~Eq. \eqref{eqn:defR}, we obtain
		\begin{align*}
			a_L \mathpzc{Z}^{(\omega)}(x)- \mathpzc{R}^{(\omega)}(x) &= a_L \Bigl(1 + \int_{0}^{x}\mathbb{W}(x; y)\Xi_\phi(y)^{-1}\omega(y)\mathpzc{Z}^{(\omega)}(y)\mathrm{d}y \Bigr)\\
			&\quad- \int_{0}^{x}\mathbb{W}(x; y)\Xi_\phi(y)^{-1}f(y)\mathrm{d}y
			- \int_{0}^{x}\mathbb{W}(x; y)\Xi_\phi(y)^{-1} \omega(y)\mathpzc{R}^{(\omega)}(x)\mathrm{d}y.
		\end{align*}
	Since the form of the above equation is the same with~Eq. \eqref{eqn:L},
	we conclude that
		\begin{equation*}
			\mathpzc{L}^{(\omega)}f(x) = a_L \mathpzc{Z}^{(\omega)}(x)- \mathpzc{R}^{(\omega)}(x).
			\label{eqn:LRelation}
		\end{equation*}
Finally, employing the 	with boundary condition $\mathpzc{L}^{(\omega)}f(a) = 0$, we get that
		\begin{equation*}
			\mathpzc{L}^{(\omega)}f(x) = \int_{0}^{a}f(y) \Xi_\phi(y)^{-1}
			\Bigl(\frac{\mathpzc{Z}^{(\omega)}(x)}{\mathpzc{Z}^{(\omega)}(a)}\mathpzc{W}^{(\omega)}(a,y) - \mathpzc{W}^{(\omega)}(x,y)\Bigr)\mathrm{d}y,
			\label{eqn:Lfx}
		\end{equation*}
	which completes  the proof.
\end{proof}

\appendix
\section{Appendix}
\subsection{Multi-refracted reflected L\'evy process}\label{multi-reftract}
In this section we derive results for a reflected multi-refracted SNLP that are used in Section \ref{reflected}. Within this subsection, we consider   $\phi(x)\equiv\phi_k(x) = \sum_{i=1}^{k}\delta_i \mathbf{1}_{(x>b_i)}$, $k\geq 1$, with $0 < \delta_1, \ldots, \delta_k$, and $-\infty:= b_0 < b_1 < \cdots < b_k < b_{k+1}:= \infty$.  We define the multi-refracted reflected Lévy process as follows. While the process is above the level $b_k$ a linear drift at rate $\delta_k$ is subtracted from the increments of process. On the other hand, when it attempts to down-cross 0, it is pushed upward so that it will not go below 0.

Under $\mathbb{P}_x$, the process, $V_k$, can be formally constructed by the recursive algorithm given below:
\begin{algorithmic}
	\item[$\textbf {Step 0}:$] For  fixed $k\geq 0$,
	set $V_{k}(0)=x$. If 
\(0< x<b_1\) or  
$b_{i} \leq x<b_{i+1}$ for  $i = 1, 2, \ldots, k$,
    then set 
    $\underline{  \mathpzc T}_0:=0$ or 
    $\underline{  \mathpzc T}_i:=0$, 
     respectively,
      and go to  \textbf {Step 1}.
	Otherwise, set $\overline {  \mathpzc   T}_0:=0$ and go to \textbf{Step 2}.
	\item[$\textbf {Step 1}:$] Let $\{\widetilde{A}_k(t) ; \; t \geq \underline{  \mathpzc   T}_i\}$ be the multi-refracted Lévy process that starts at the time $\underline{  \mathpzc   T}_i$ at the level $x$ and it is a solution of~\eqref{eqn:levelDependent} with $\phi = \phi_k $ (which exists by Theorem 1~\cite{CPRY2019}). 	
	Three different possibilities may occur at the first exit time of $\widetilde{A}_{k}(t)$ from
	$[b_i, b_{i+1}]$. \\
	\begin{itemize}
		\item [(i)] Let $\underline{  \mathpzc T}_{i+1}=\inf\{t> \underline{  \mathpzc T}_{i}:\widetilde{A}_{k}(t)>b_{i+1} \}$. Set $V_k(t)=\widetilde A_k(t)$ for  $\underline{  \mathpzc   T}_i\leq t< \underline{ \mathpzc   T}_{i+1}$ and then go to \textbf{Step 1} with $i+1$.
		\item [(ii)] Let $\widehat{\mathpzc T}_{i}=\inf\{t>\underline{\mathpzc{T}}_i: 0<\widetilde A_k(t)<b_i\}$ and  \( n_{i,k} = \min\{j\leq i: \widetilde{A}_k(\widehat{\mathpzc T}_{i})> b_{i-j}\mathbf{1}_{(i>j)}\} \). Set
		$V_k(t)=\widetilde A_k(t)$ for  $\underline {\mathpzc{T}}_i\leq t<\widehat{\mathpzc T}_{i}$ and then go to \textbf{Step 1}  replace solely  \(\underline{\mathpzc T}_i \) with \(\widehat{\mathpzc T}_{i}\) and take the index $i-n_{i,k}$ in the rest of the terms.
		\item [(iii)] Let $\overline {\mathpzc T}_i=\inf\{t>\underline {\mathpzc T}_i:\widetilde A_k(t)<0\}$. Set  $V_k(t)=\widetilde A_k(t)$  for $\underline{\mathpzc T}_i\leq t <\overline {\mathpzc T}_i$ and then go to \textbf{Step 2} (with $i$).
	\end{itemize}
	\item[$\textbf {Step 2}:$] Let $\{\widetilde{U}(t) ; \; t \geq \overline{  \mathpzc   T}_i\}$ be the L\'evy process reflected at the lower boundary 0 that starts at time $\overline {\mathpzc T}_i$ at 0.  Let  $\widetilde {\mathpzc T}_i=\inf\{t>\overline{\mathpzc T}_i:\widetilde U(t)>b_1\}$. Set $V_k(t)=\widetilde U(t)$  for $\overline{\mathpzc T}_i\leq t <\widetilde {\mathpzc T}_i$ and $x=b_1$. Then,  go to \textbf{Step 1}  and replace solely  $\underline {\mathpzc T}_i $ with $\widetilde {\mathpzc T}_i$ and start with the index $i=1$ in the rest of the terms.
\end{algorithmic}
Let us now define the non-decreasing and right-continuous process $R_k(t)$ that represents the cumulative amounts of modification up to $t$ that pushes the process upward when it attempts to go below 0,  from level $k$. Then we have the decomposition
	\begin{equation}
		\mathrm{d}V_k(t) = \mathrm{d}X(t) + \mathrm{d}R_k(t)
		-\sum_{i=1}^{k}\delta_i\mathbf{1}_{(V_k(t)>b_i)}\mathrm{d}t.
		\label{eqn:RefMulti}
	\end{equation}
In particular, in the case of bounded variation
\[R_k(t)=\sum_{t\geq 0,V_k(t-)+\Delta X_t<0}\big|V_k(t-)+\Delta X_t\big|, \quad t\geq 0,\]
where $\Delta \xi_t=\xi_t-\xi_{t-}$, \(t\geq0\),  for any right continuous process \(\xi\).

   This extends the case of \cite{PY2018, PY2017} in a multi-refracted framework and   the  multi-refracted case of \cite{CPRY2019} in a reflected environment.

We start showing that for $V_k(t)$   a strong approximation
exists between the bounded and the unbounded variation case, which  guarantees
the existence
of the solution in Eq. \eqref{eqn:RefMulti}. To do this,  we shall  use \textbf{(i)} the fact that [see pp. 210 in~\cite{Bertoin1996}] for any spectrally negative L\'evy process with unbounded variation
paths, $X$, there exists a sequence of bounded variation $X^{(n)}$ such that
\[\lim_{n\to \infty}\sup_{0 \leq s \leq t}|X^{(n)}(s) - X(s)| = 0 \;  \mbox{ a.s. for each } \; t >0,\]
and \textbf {(ii)} the lemma below,  (which for $k=1$ is equivalent to the  Lemma A.1 in \cite{PY2018} and thus the proof is omitted).
\begin{lemma}	\label{lem:A1}
	For $k \geq 1$ and fixed $t > 0$, let $(x_s)_{0 \leq s \leq t}$ and
	$(\tilde{x}_s)_{0 \leq s \leq t}$ be the paths of two
	different L\'evy processes such that
$
		\sup_{0\leq s\leq t}|x_s - \tilde{x}_s| < \varepsilon, \quad \mbox{for some } \varepsilon > 0.$
	Also, fix $z_k, \tilde{z}_k \in \R$ and $0 \leq t_0 < t$.
	\begin{itemize}
		\item[\upshape(i)] Define the reflected paths $y_s(z_k, t_0)$ and $\tilde{y}_s(\tilde{z}_k, t_0)$
		on $[t_0, t]$ of the shifted paths
$ z_k + (x_s - x_{t_0})$ and $  \tilde{z}_k +(\tilde{x}_s - \tilde{x}_{t_0}),$
		respectively,  i.e.  for all $t_0 \leq s \leq t$, let
			\begin{align*}
				& y_s(z_k, t_0):= z_k+(x_s-x_{t_0})+\bigl(-\inf _{t_0 \leq u \leq s}\bigl[z_k+\left(x_u-x_{t_0}\right)\bigr]\bigr) \vee 0, \\
				& \tilde{y}_s(\tilde{z}_k, t_0) :=\tilde{z}_k+ (\tilde{x}_s-\tilde{x}_{t_0})+ \bigl(-\inf _{t_0 \leq u \leq s}\bigl[\tilde{z}_k+\left(\tilde{x}_u-\tilde{x}_{t_0}\right)\bigr]\bigr) \vee 0.
			\end{align*}
		Then, we have
			\begin{equation*}
				\sup _{t_0 \leq s \leq t}\left|y_s(z_k, t_0)-\tilde{y}_s(\tilde{z}_k, t_0)\right|<2|z_k-\tilde{z}_k|+4 \varepsilon .
			\end{equation*}
		\item[\upshape(ii)] Similarly, for the multi-refracted paths $a_s(z_{k-1}, t_0)$ and $\tilde{a}_s(\tilde{z}_{k-1}, t_0)$
		for all $t_0 \leq s \leq t$ that solves
			\begin{align*}
				& a_s(z_{k-1}, t_0):= z_{k-1} +(x_s-x_{t_0})-\sum_{i=1}^{k}\delta_i \int_{t_0}^{s} \mathbf 1_{\{a_u(z_{k-1}, t_0)>b_i\}} \mathrm{d} u, \\
				& \tilde{a}_s(\tilde{z}_{k-1}, t_0) := \tilde{z}_{k-1} +(\tilde{x}_s-\tilde{x}_{t_0})-\sum_{i=1}^{k}\delta_i \int_{t_0}^{s} \mathbf 1_{\{\tilde{a}_u(\tilde{z}_{k-1}, t_0)>b_i\}} \mathrm{d}u,
			\end{align*}
		it holds that
			\begin{equation*}
				\sup _{t_0 \leq s \leq t}\left|a_s(z_{k-1}, t_0)-\tilde{a}_s(\tilde{z}_{k-1}, t_0)\right|<2|z_{k-1}-\tilde{z}_{k-1}|+4 \varepsilon .
			\end{equation*}
	\end{itemize}
\end{lemma}
\begin{proposition}	\label{prop:ConvergenceVk}
	Assume that $X$ is of unbounded variation and $(X^{(n)})_{n\geq 1}$ is a strongly approximating sequence of bounded variation. In addition, we define $V_k$ and $V_k^{(n)}$ as the reflected
	multi-refracted processes associated with $X$ and $X^{(n)}$, respectively. Then, $V_k^{(n)}$
	is the strongly approximating
	sequence of $V_k$.
\end{proposition}
\begin{proof}
	The proof extends  the line of logic as of  \cite{PY2018} for the process in  \eqref{eqn:RefMulti}. It is essential to introduce and work with a  new   broader family of  stopping times  in order to separate  the refracted  from the reflected paths.
		Let $\beta_k:= b_k/2 > 0$ for $k \geq 1$ and  define a sequence of increasing random times,   $\overline{\tau}_{0}:= 0$ (for convenience) and 	$(\underline{\tau}_{1,k}$, $\overline{\tau}_1$, $\underline{\tau}_{2,k}$, $\overline{\tau}_2$, $\ldots$, $\underline{\tau}_{\nu,k}$, $\overline{\tau}_\nu )$ as follows
		\begin{align*}
			\underline{\tau}_{1,k} &:= \inf\{s>0: V_{k-1}(s) > \beta_k\}, \quad k = 1, 2, \ldots, \\
			\sigma_{1, k} &:= \inf\{s>\underline{\tau}_{1,k}: V_{k}(s) = 0 \} \; \mbox{with } \; \sigma_{1, 0} := \inf\{s> \underline{\tau}_{1,1}: V_{k}(s) = 0 \}, \quad k = 0, 1, 2,\ldots,\\
			\widetilde{\sigma}_{1} &:= \inf\{k \geq 0: \sigma_{1, k}\} = \inf\{\sigma_{1, 0},
			\sigma_{1, 1}, \ldots, \sigma_{1, k}\},\\
			\overline{\tau}_{1,k} &:= \sup\{s <\widetilde{\sigma}_{1}: V_{k-1}(s) > \beta_k \}, \quad k = 1, 2, \ldots, \\
			\overline{\tau}_{1} &:= \sup\{k \geq 0: \overline{\tau}_{1, k}\} = \sup\{\overline{\tau}_{1, 1},
			\overline{\tau}_{1, 2}, \ldots, \overline{\tau}_{1, k}\},
		\end{align*}
	and for all $\nu \geq 2$
		\begin{align*}
			\underline{\tau}_{\nu,k} &:= \inf\{s>:\widetilde{\sigma}_{\nu-1}: V_{k-1}(s) > \beta_k\}, \quad k = 1, 2, \ldots, \\
			\sigma_{\nu, k} &:= \inf\{s>\underline{\tau}_{\nu,k}: V_{k}(s) = 0 \} \;
			\mbox{with } \; \sigma_{\nu, 0} := \inf\{s> \underline{\tau}_{\nu,1}: V_{k}(s) = 0 \}, \quad k = 0, 1, 2,\ldots,\\
			\widetilde{\sigma}_{\nu} &:= \inf\{k \geq 0: \sigma_{\nu, k}\} = \inf\{\sigma_{\nu, 0},
			\sigma_{\nu, 1}, \ldots, \sigma_{\nu, k}\},\\
			\overline{\tau}_{{\color{blue}\nu},k} &:= \sup\{s <\widetilde{\sigma}_{\nu}: V_{k-1}(s) > \beta_k \}, \quad k = 1, 2, \ldots, \\
			\overline{\tau}_{\nu} &:= \sup\{k \geq 0: \overline{\tau}_{\nu, k}\} = \sup\{\overline{\tau}_{\nu, 1},
			\overline{\tau}_{\nu, 2}, \ldots, \overline{\tau}_{\nu, k}\},
		\end{align*}
	Further, let
	\(N = 1 + N_1 + N_2, \; \text{with}\; N_1 := \sup\{ \nu \geq 0: \underline{\tau}_{\nu,1} < t\} \; \mbox{and } \;
	N_2 := \sup\{ \nu \geq 0: \overline{\tau}_{\nu} < t\}\)
 be the total number of times switching has occurred until time $t$ (plus one ) and define
		\begin{equation*}
			\underline{\beta}_k := \min_{1 \leq \nu \leq N_1} \inf_{s \in [\underline{\tau}_{\nu,1},\overline{\tau}_{\nu}\wedge t )}V_k(s), \;
			\;\;  \;
			\underline{\beta} = \min\{k\geq 0: \underline{\beta}_k \},
		\end{equation*}
	where $\underline{\beta} > 0$, by the definition of $\overline{\tau}$ and $\underline{\tau}$.
	
It suffices to show that  there exists a finite $C$ and
	$\underline{n} \in \N$ such that
		\begin{equation}
			\sup _{0 \leq s \leq t}|V_k^{(n)}(s)-V_k(s)| \leq C \sup _{0 \leq s \leq t}|X^{(n)}(s)-X(s)|, \quad n \geq \underline{n}.
			\label{eqn:strongAp1}
		\end{equation}
We shall  choose $\underline{n}$ large enough so that

		\begin{equation} \label{eq2}
			[4(2^{N}-1)]\sup_{m \geq \underline{n}}\sup_{0 \leq s \leq t} |X^{(m)}(s)-X(s)|<  \underline{\beta},
		\end{equation}
	which we will see in later that, for all $n \geq \underline{n}$, this bound
	confirms that $\underline{\tau}_{\nu,1}$ and $\overline{\tau}_{\nu}$
	can act as switching times for both $V_k(t)$ and $V_k^{(n)}(t)$ i.e. on each interval
	$[\underline{\tau}_{\nu,1}, \overline{\tau}_{\nu})$
	and $[\overline{\tau}_{\nu}, \underline{\tau}_{\nu+1,1})$ are  multi-refracted paths
	and reflected paths for both $V_k(t)$ and $V_k^{(n)}(t)$, respectively.
	
	Let us fix $n > \underline{n}$ and $\varepsilon := \sup _{0 \leq s \leq t}|X^{(n)}(s)-X(s)|$ and define  a sequence $(\eta_{\nu})_{0 \leq  \nu \leq N}$ such that $\eta_0 = 0$ and  \( \eta_{\nu+1} = 2 \eta_{\nu} + 4 \varepsilon\), which gives
\(		 \eta_{\nu} = 4(2^{N}-1)\varepsilon\),
	and by Eq. \eqref{eq2}, we have
		\[4(2^{N}-1)\varepsilon< \underline{\beta} < b_1/2. \]
	Then, we shall show that the followings for reflected and multi-refracted paths of the process
		\begin{align}
			\sup _{\overline{\tau}_{\nu} \leq s \leq  \underline{\tau}_{\nu+1,1} \wedge t} &|V_k(s)-V_k^{(n)}(s)| \leq \eta_{2{\color{blue}\nu} +1}, \; \mbox{a.s. } \; 0\leq  \nu \leq N_2,
			\label{eqn:reflectedPaths}\\
			\sup _{\underline{\tau}_{\nu,1} \leq s \leq \overline{\tau}_{\nu} \wedge t} &|V_k(s)-V_k^{(n)}(s)| \leq \eta_{2{\color{blue}\nu}}, \; \mbox{a.s. }\; 0\leq  \nu \leq N_1,
			\label{eqn:multirefracPaths}
		\end{align}
	and hence~\eqref{eqn:strongAp1} holds true with $C = 4(2^{N}-1)$. Thus, the results above lead
	us to prove the following claims.
	
	\noindent \textbf{Claim 1:}
	Fix $\nu \geq 0$.
	Assume that $\overline{\tau}_{\nu} < t$ (or $\nu \leq N_2$),  $\tilde{\eta} := 2 \eta + 4 \varepsilon < \underline{\beta} < b_1/2$ and
	\[ |V_k(\overline{\tau}_{\nu^-})+ \Delta X(\overline{\tau}_{\nu})- (V_k^{(n)}(\overline{\tau}_{\nu^-}) + \Delta X^{(n)}(\overline{\tau}_{\nu}))| < \eta. \]
	Then, we have
		\[|V_k(s)-V_k^{(n)}(s)| \leq \tilde{\eta}, \quad \mbox{for } \quad \overline{\tau}_{\nu} \leq s \leq  \underline{\tau}_{\nu+1,1} \wedge t.\]
	\begin{proof}
		Consider the reflected paths on $[\overline{\tau}_{\nu}, \underline{\tau}_{\nu+1,1} \wedge t]$
			\begin{align*}
				Y_k(s) &:= V_k(\overline{\tau}_{\nu^-}) + \Delta X(\overline{\tau}_{\nu}) + (X(s)- X(\overline{\tau}_{\nu}))+ \bigl[ - \inf_{ \overline{\tau}_{\nu} \leq u \leq s} \bigl(V_k(\overline{\tau}_{\nu^-}) + \Delta X(\overline{\tau}_{\nu}) + (X(s)- X(\overline{\tau}_{\nu})) \bigr)\bigr]\vee 0,
				 			\end{align*}
			 			and
			 						\begin{align*}
				Y_k^{(n)}(s) &:= V_k^{(n)}(\overline{\tau}_{\nu^-}) + \Delta X^{(n)}(\overline{\tau}_{\nu}) + (X^{(n)}(s)- X^{(n)}(\overline{\tau}_{\nu})) \\
				&\qquad+ \bigl[ - \inf_{ \overline{\tau}_{\nu} \leq u \leq s} \bigl(V_k^{(n)}(\overline{\tau}_{\nu^-}) + \Delta X^{(n)}(\overline{\tau}_{\nu}) + (X^{(n)}(s)- X^{(n)}(\overline{\tau}_{\nu})) \bigr)\bigr]\vee 0.
			\end{align*}
		Applying Lemma \ref{lem:A1} (i) for $z_k =  V_k(\overline{\tau}_{\nu^-}) + \Delta X(\overline{\tau}_{\nu})$ and $\tilde{z}_k =  V_k^{(n)}(\overline{\tau}_{\nu^-}) + \Delta X^{(n)}(\overline{\tau}_{\nu})$ and $t_0 = \overline{\tau}_{\nu}$,
		we obtain that
			\begin{align*}
				|Y_k(s) - Y_k^{(n)}(s)| &< 2 \big|V_k(\overline{\tau}_{\nu^-}) + \Delta X(\overline{\tau}_{\nu})- (V_k^{(n)}(\overline{\tau}_{\nu^-})+ \Delta X^{(n)}(\overline{\tau}_{\nu}))\big|+ 4 \varepsilon\leq 2 \eta + 4 \varepsilon = \tilde{\eta} < \underline{\beta},
			\end{align*}
		for all $s \in [\overline{\tau}_{\nu}, \underline{\tau}_{\nu+1,1} \wedge t]$. Since $V_k\leq \beta$ for $[\overline{\tau}_{\nu}, \underline{\tau}_{\nu+1,1}]$ and using the fact that \(\tilde{\eta} := 2 \eta + 4 \varepsilon < \underline{\beta} < b_1/2\) we can conclude that there is no refraction on $[\overline{\tau}_{\nu}, \underline{\tau}_{\nu+1,1}\wedge t]$.
 Therefore,
		$V_k(s)$ and $V_k^{(n)}(s)$ coincide with defined $Y_k(s)$ and $Y_k^{(n)}(s)$ above on  $[\overline{\tau}_{\nu}, \underline{\tau}_{\nu+1,1}\wedge t]$ and thus the claim is proved.
	\end{proof}
	
\noindent 	\textbf{Claim 2:}
	Fix $\nu \geq 0$ and $k \geq 1$.
	Assume that $\underline{\tau}_{\nu,1} < t$ (or $\nu \leq N_2$), $\tilde{\eta} := 2 \eta + 4 \varepsilon < \underline{\beta} < b_1/2$ and
		\begin{equation*}
			|V_{k-1}(\underline{\tau}_{\nu,k})-V_{k-1}^{(n)}(\underline{\tau}_{\nu,k})| \leq \eta.
		\end{equation*}
	Then,
	\begin{align*}
		|V_k(s)-V_k^{(n)}(s)| \leq \tilde{\eta}, \quad \mbox{for} \quad  \underline{\tau}_{\nu,1}
		\leq s \leq  \overline{\tau}_{\nu} \wedge t,\\
		|V_k(\overline{\tau}_{\nu^-})+ \Delta X(\overline{\tau}_{\nu})- (V_k^{(n)}(\overline{\tau}_{\nu^-}) +\Delta X^{(n)}(\overline{\tau}_{\nu}))| < \tilde{\eta}, \quad \mbox{if } \; \overline{\tau}_{\nu}< t.
		\end{align*}
The proof of \textbf{Claim 2} follows similar line of logic as in \textbf{Claim 1} and use of Lemma \ref{lem:A1} (ii), and thus is omitted. 	
	
	Now, we will show Eqs.~(\ref{eqn:reflectedPaths})-(\ref{eqn:multirefracPaths}) by using mathematical induction. We start with $\nu =0$ by Claim 1, it is obvious that
		\(|V_k(0^-) + \Delta X(0)- (V_k(0^-) + \Delta X(0))| = |x-x| = \eta_0 = 0,\)
	since for convenience we take $\overline{\tau}_{0} := 0$.
	Then, applying repeatedly Claim 1 and Claim 2 one after other for $N$ times, Eqs.~(\ref{eqn:reflectedPaths})-(\ref{eqn:multirefracPaths}) hold true.
	
	Finally, repeating the same reasoning given in  Lemma 4.1 in \cite{PY2018}, we verify that for all driving processes $X$ of
	unbounded variation when $x$ is fixed we have, for $1 \leq j \leq k$
		\begin{equation*}
			\Prob_x(V_k^{(\infty)}(t)= b_j) = 0 \; \mbox{for almost every } t \geq 0,
			\label{eqn:diffuseMeasureVk}
		\end{equation*}
	which leads us to obtain, with \(R_k^{(n)}(t)\) as before for the process \(V_k^{(n)}\), that
		\begin{align*}
			\lim_{n\to \infty}V_k^{(n)}(t) &= \lim_{n \to \infty}\Bigl(X^{(n)}(t) + R_k^{(n)}(t)
			-\sum_{i=1}^{k}\delta_i\int_{0}^{t}\mathbf{1}_{(V_k^{(n)}(s)>b_i)}\mathrm{d}s\Bigr)= X_t + R_k^{(\infty)}
			-\sum_{i=1}^{k}\delta_i\int_{0}^{t}\mathbf{1}_{(V_k^{(\infty)}(s)>b_i)}\mathrm{d}s.
		\end{align*}
 This result gives
	$V_k := V_k^{(\infty)}$ solves Eq. \eqref{eqn:RefMulti}, which completes the proof.
\end{proof}
Next, we provide identities for the fluctuations of  $V_k$ defined in Eq.  \eqref{eqn:RefMulti}. To do this, we denote by $\mathbb{W}_k^{(q)}$, $\mathbb{Z}_k^{(q)}$ the corresponding scale functions of Eqs  \eqref{eqn:shiftedLevelScaleW}-\eqref{eqn:shiftedLevelScaleZ}  for $\phi(x)\equiv\phi_k(x) = \sum_{j=1}^{k}\delta_j \mathbf{1}_{(x>b_j)}$, for which the following identities are true.
\begin{remark}
	\label{rem:remark34}	
	For $0 \leq x \leq a$ and $d < b_1 < \cdots <  b_k  \leq  a$ and $\mathbb{Z}_{k}^{(q)}(x)\equiv\mathbb{Z}_{k}^{(q)}(x,0)$ it holds that
		\begin{align*}
			\int_0^{a-b_k} \int_{(y,\infty)} &\mathbb{W}_{k-1}^{(q)}(b_k+y-\theta ; d) \Bigl[\frac{W_k^{(q)}\left(x-b_k\right) W_k^{(q)}(a-b_k-y)}{W_k^{(q)}(a-b_k)}-W_k^{(q)}(x-b_k-y)\Bigr] \nonumber  \nu(\mathrm{d} \theta) \mathrm{d}y\nonumber\\
			&=-\frac{W_k^{(q)}(x-b_k)}{W_k^{(q)}(a-b_k)}\Bigl(\mathbb{W}_{k-1}^{(q)}(a ; d)+\delta_k \int_{b_k}^a W_k^{(q)}(a-y) \mathbb{W}_{k-1}^{(q){\prime}}(y ; d) \mathrm{d} y\Bigr) \nonumber\\
			&\quad +\mathbb{W}_{k-1}^{(q)}(x ; d)+\delta_k \int_{b_k}^x W_k^{(q)}(x-y) \mathbb{W}_{k-1}^{(q){\prime}}(y ; d) \mathrm{d}y.
			\label{eqn:OUWk}
		\end{align*}
	Similarly, we have
		\begin{align*}
		&	\int_0^{a-b_k} \int_{y}^{\infty} \mathbb{Z}_{k-1}^{(q)}(b_k+y-\theta) \Bigl[\frac{W_k^{(q)}\left(x-b_k\right) W_k^{(q)}\left(a-b_k-y\right)}{W_k^{(q)}\left(a-b_k\right)}-W_k^{(q)}(x-b_k-y)\Bigr]
			  \nu(\mathrm{d} \theta) \mathrm{d}y\nonumber\\
			& =-\frac{W_k^{(q)}(x-b_k)}{W_k^{(q)}(a-b_k)}\Bigl(\mathbb{Z}_{k-1}^{(q)}(a)+\delta_k \int_{b_k}^a W_k^{(q)}(a-y) \mathbb{Z}_{k-1}^{(q){\prime}}(y)
			\mathrm{d} y\Bigr)+\mathbb{Z}_{k-1}^{(q)}(x)+\delta_k \int_{b_k}^x W_k^{(q)}(x-y) \mathbb{Z}_{k-1}^{(q){\prime}}(y) \mathrm{d}y.
		\end{align*}
\end{remark}
\noindent The first eq. of Remark \ref{rem:remark34} can be found in Remark 34 in \cite{CPRY2019}. The second eq.  can be proved by employing Eq. \eqref{eqn:Zinduction} (see below) and perform algebraic manipulations.

\noindent Let
	\begin{equation*}
		T_{k}^{a, +} = \inf\{t > 0: V_k(t) > a\}, \quad \mbox{and} \quad
		T_{k}^{a, -}  = \inf\{t > 0: V_k(t) < a\},
		\label{eqn:passageTimesMultiRef}
	\end{equation*}
be the first passage times for the reflected multi-refracted process for $a > 0$ and $k \geq 1$.
\begin{theorem}\label{thm:refMultiResolvent}
	For $q\geq 0$, $y \leq x \leq a$, and $y < b_1 < b_2 < \cdots < b_k  \leq  a$, and a Borel set $\mathfrak{B}$ on $[0, a]$,
	we have
		\begin{equation}\label{eqn:refMultiResolvent}
			\Ex_x\Bigl[\int_{0}^{T_{k}^{a, +} }e^{-qt}\mathbf{1}_{(V_k(t) \in \mathfrak{B})}\mathrm{d}t\Bigr]
			=\int_{\mathfrak{B}} \Xi_{\phi_k}(y)^{-1}\Bigl(\mathbb{W}_k^{(q)}(a; y)\frac{\mathbb{Z}_k^{(q)}(x)}{\mathbb{Z}_k^{(q)}(a)}-\mathbb{W}_k^{(q)}(x; y)\Bigr)\mathrm{d}y.
		\end{equation}
\end{theorem}
\begin{proof}
	\textbf{(i) Bounded Variation:}  The proof is inductive. Let  $f^{(q)}(x,a, \mathfrak{B})$ be expectation of the left hand side in the above theorem.  For $k=1$,
	 Eq. \eqref{eqn:refMultiResolvent} agrees with Theorem 4.1 in \cite{PY2018} (after using Eq. (32) of \cite{CPRY2019}).
	Then,  using the induction step and the  strong Markov property we have for $x <  b_k$ that
		\begin{align}
			f^{(q)}(x,a, \mathfrak{B}) &= \Ex_x\Bigl[\int_{0}^{T_{k-1}^{b_k, +} }e^{-qt}\mathbf{1}_{(V_{k-1}(t) \in \mathfrak{B})}\mathrm{d}t\Bigr] +
			\Ex_x\Bigl[e^{-qT_{k-1}^{b_k, +}}\mathbf{1}_{(T_{k-1}^{b_k, +} < \infty)}\Bigr] f^{(q)}(b_k,a, \mathfrak{B})\nonumber\\
			&=\int_{0}^{\infty}\sum_{i=0}^{k-1}\frac{\mathbb{W}_{k-1}^{(q)}(b_k; y)\frac{\mathbb{Z}_{k-1}^{(q)}(x)}{\mathbb{Z}_{k-1}^{(q)}(b_k)}-\mathbb{W}_{k-1}^{(q)}(x; y)}{\Xi_{\phi_i}(y)} \mathbf{1}_{(y \in (b_i, b_{i+1}])}\mathrm{d}y + \frac{\mathbb{Z}_{k-1}^{(q)}(x)}{\mathbb{Z}_{k-1}^{(q)}(b_k)}f^{(q)}(b_k,a, \mathfrak{B}).
			\label{eqn:refMulti1}
		\end{align}
	where the last eq. holds after using 	\(
	\Ex_{x}[e^{-qT_{k-1}^{b_k, +}}\mathbf{1}_{(T_{k-1}^{b_k, +} < \infty)}]={\mathbb{Z}_{k-1}^{(q)}(x)}/{\mathbb{Z}_{k-1}^{(q)}(b_k)},
	\)
 which is true by inductive hypothesis.
	
\noindent 	On the other hand, for $b_k \leq x \leq a$, again  using strong Markov property  gives
		\begin{align*}
			f^{(q)}(x,a, \mathfrak{B}) &= \Ex_x\Bigl[\int_{0}^{\tau_{k}^{a, +} \wedge \tau_{k}^{b_k, -}  }e^{-qt}\mathbf{1}_{(X_{k}(t) \in \mathfrak{B})}\mathrm{d}t\Bigr]+\Ex_x\Bigl[e^{-q\tau_{k}^{b_k, -}}f^{(q)}(X_k(\tau_{k}^{b_k, -}), a, \mathfrak{B})\mathbf{1}_{(\tau_{k}^{b_k, -} < \tau_{k}^{a, -})}\Bigr]\\
			&= \int_{0}^{\infty}\Bigl[\frac{W_k^{(q)}(x-b_k) W_k^{(q)}(a-y)}{W_k^{(q)}(a-b_k)}-W_k^{(q)}(x-y)\Bigr]
			\mathbf{1}_{(y \in (b_k,a])}\mathrm{d}y\\
				&\quad+ \int_{0}^{\infty}\int_{z}^{\infty}f^{(q)}(z-\theta+b_k, a, \mathfrak{B})
			\Bigl[\frac{W_k^{(q)}(x-b_k) W_k^{(q)}(a-b_k-z)}{W_k^{(q)}(a-b_k)}-W_k^{(q)}(x-b_k-z)\Bigr]\nu(\mathrm{d}\theta)\mathrm{d}z.
		\end{align*}
	Then, substituting Eq. \eqref{eqn:refMulti1} in the  second term of  the above equation, we have
		\begin{align*}
			f^{(q)}(x,a, \mathfrak{B}) &= \int_{0}^{\infty}\Bigl[\frac{W_k^{(q)}(x-b_k) W_k^{(q)}(a-y)}{W_k^{(q)}(a-b_k)}-W_k^{(q)}(x-y)\Bigr]
			\mathbf{1}_{(y \in (b_k,a])}\mathrm{d}y\\
			&\quad+\int_{0}^{\infty} \int_{z}^{\infty}\int_{0}^{\infty}\frac{\mathbb{W}_{k-1}^{(q)}(b_k; y)}{\mathbb{Z}_{k-1}^{(q)}(b_k)}\sum_{i=0}^{k-1}
			\frac{\mathbb{Z}_{k-1}^{(q)}(z-\theta+b_k)}{\Xi_{\phi_i}(y)}\mathbf{1}_{(y \in (b_i, b_{i+1}])}\mathrm{d}y\\
			&\qquad \times  \Bigl[\frac{W_k^{(q)}(x-b_k) W_k^{(q)}(a-b_k-z)}{W_k^{(q)}(a-b_k)}-W_k^{(q)}(x-b_k-z)\Bigr]
			\nu(\mathrm{d}\theta)\mathrm{d}z\\
			&\quad- \int_{0}^{\infty}\int_{z}^{\infty}\int_{0}^{\infty}\sum_{i=0}^{k-1}\frac{1}{\Xi_{\phi_i}(y)}\mathbb{W}_{k-1}^{(q)}(z-\theta+b_k; y)\mathbf{1}_{(y \in (b_i, b_{i+1}])}\mathrm{d}y\\
			&\qquad \times \Bigl[\frac{W_k^{(q)}(x-b_k) W_k^{(q)}(a-b_k-z)}{W_k^{(q)}(a-b_k)}-W_k^{(q)}(x-b_k-z)\Bigr]
			\nu(\mathrm{d}\theta)\mathrm{d}z\\		
			&\quad + \frac{f^{(q)}(b_k,a,\mathfrak{B})}{\mathbb{Z}_{k-1}^{(q)}(b_k)}
			\int_{0}^{\infty}\int_{z}^{\infty}\mathbb{Z}_{k-1}^{(q)}(z-\theta +b_k)\\
			&\qquad \times \Bigl[\frac{W_k^{(q)}(x-b_k) W_k^{(q)}(a-b_k-z)}{W_k^{(q)}(a-b_k)}-W_k^{(q)}(x-b_k-z)\Bigr]\nu(\mathrm{d}\theta)\mathrm{d}z.
		\end{align*}
	Next, using Remark \ref{rem:remark34}, we simplify the above expression to
		\begin{align*}
			f^{(q)}(x,a, \mathfrak{B}) &=  \int_{0}^{\infty}\Bigl[\frac{W_k^{(q)}(x-b_k) W_k^{(q)}(a-y)}{W_k^{(q)}(a-b_k)}-W_k^{(q)}(x-y)\Bigr]
			\mathbf{1}_{(y \in (b_k,a])}\mathrm{d}y\\
			&\quad +\int_{0}^{\infty}\sum_{i=0}^{k-1}\frac{1}{\Xi_{\phi_i}(y)}
			\Bigl(\frac{W_k^{(q)}(x-b_k)}{W_k^{(q)}(a-b_k)}\mathbb{W}_{k}^{(q)}(a; y) -\mathbb{W}_{k}^{(q)}(x; y)\Bigr)\mathbf{1}_{(y \in (b_i, b_{i+1}])}\mathrm{d}y\\
			&\quad +\int_{0}^{\infty}\sum_{i=0}^{k-1}\frac{1}{\Xi_{\phi_i}(y)}
			\Bigl[ \frac{\mathbb{W}_{k-1}^{(q)}(b_k; y)}{\mathbb{Z}_{k-1}^{(q)}(b_k)}
			\Bigl(-\frac{W_k^{(q)}(x-b_k)}{W_k^{(q)}(a-b_k)}\mathbb{Z}_{k}^{(q)}(a) -\mathbb{Z}_{k}^{(q)}(x)\Bigr)\Bigr]\mathbf{1}_{(y \in (b_i, b_{i+1}])}\mathrm{d}y\\
			&\quad+\frac{f^{(q)}(b_k,a,\mathfrak{B})}{\mathbb{Z}_{k-1}^{(q)}(b_k)}
			\Bigl(-\frac{W_k^{(q)}(x-b_k)}{W_k^{(q)}(a-b_k)}\mathbb{Z}_{k}^{(q)}(a) -\mathbb{Z}_{k}^{(q)}(x)\Bigr).
		\end{align*}
Finally,  setting $x = b_k$ and noticing that $\mathbb{W}_{k-1}^{(q)}(b_k) = \mathbb{W}_{k}^{(q)}(b_k)$ and $\mathbb{Z}_{k-1}^{(q)}(b_k) = \mathbb{Z}_{k}^{(q)}(b_k)$ leads to a explicit formula for $f^{(q)}(b_k, a, \mathfrak{B})$.
	Finally, putting the expression for $f^{(q)}(b_k, a, \mathfrak{B})$ into Eq. \eqref{eqn:refMulti1} and using the fact that  \(\mathbb{W}_k^{(q)}(x; y)=\Xi_{\phi_k}(y) W_k^{(q)}(x-y)\) for $y \in (b_k, a]$ [see Lemma 14 in \cite{CPRY2019}]
	one gets
		\begin{align*}
			f^{(q)}(x, a, \mathfrak{B}) &= \int_{0}^{\infty} \sum_{i=0}^{k}\frac{1}{\Xi_{\phi_i}(y)}\Bigl[\frac{\mathbb{Z}_{k}^{(q)}(x)}{\mathbb{Z}_{k}^{(q)}(a)}\mathbb{W}_{k}^{(q)}(a; y) -\mathbb{W}_{k}^{(q)}(x; y)
			\Bigr]\mathbf{1}_{(y \in (b_i, b_{i+1}] \bigcap [0, a))}\mathrm{d}y\\
			&=\int_{0}^{\infty}\frac{1}{\Xi_{\phi_k}(y)}\Bigl[\frac{\mathbb{Z}_{k}^{(q)}(x)}{\mathbb{Z}_{k}^{(q)}(a)}\mathbb{W}_{k}^{(q)}(a; y) -\mathbb{W}_{k}^{(q)}(x; y)
			\Bigr]\mathbf{1}_{(y \in [0, a))}\mathrm{d}y,
		\end{align*}
	which the last eq. follows by Definition 4 in \cite{CPRY2019}.

\noindent 	\textbf{(ii) Unbounded Variation:}
	We now extended the result to the case of unbounded variation.
	Let $V_k^{(n)}(t)$ is be the reflected multi-refracted  process
	associated with $X^{(n)}$  and $T_{k}^{a, (n), +}$ is the first passage time above $a$ of $V_k^{(n)}$.
	It is clear that $T_{k}^{a, (n), +}$ and $T_{k}^{a, +}$ are both finite a.s. because
	they are bounded from above by the upcrossing times at $a$ for the reflected processes
	of drift changed process $X_k(t)$ for any $k \geq 0$ (this will be also confirmed in  Theorem \ref{thm:OneSidedMultiRefRac} below). In order to prove that
\( T_{k}^{a, (n), +} \xrightarrow{n \uparrow \infty}   T_{k}^{a, +}, \quad \mbox{a.s.,}\)
	it is enough to show that \(T_{k}^{a^-, +} =T_{k}^{a, +}\) a.s., which can be obtained  following the same arguments as
	in the proof of Theorem 4.1 in~\cite{PY2018}.
	Therefore, by dominated convergence
	theorem and~Proposition \ref{prop:ConvergenceVk}, we get that
		\[ \lim_{n \to \infty}\Ex_x\Bigl[\int_{0}^{T_{k}^{a,(n), +} }e^{-qt}\mathbf{1}_{(V_k^{(n)}(t) \in \mathfrak{B})}\mathrm{d}t\Bigr]  = \Ex_x\Bigl[\int_{0}^{T_{k}^{a, +} }e^{-qt}\mathbf{1}_{(V_k(t) \in \mathfrak{B})}\mathrm{d}t\Bigr], \quad  \text{a.s.}\]
	Moreover, since  \(\Prob_x(V_k(t) \in \partial \mathfrak{B})=\Prob_x(\sup_{0 \leq s \leq t}V_k(s)= a) = 0 \) for Lebesgue a.e. $t > 0$,
	we conclude the convergence of LHS of Eq. \eqref{eqn:refMultiResolvent}.
	It remains to show the convergence of the RHS of Eq. \eqref{eqn:refMultiResolvent}. To do this recall from  Theorem 7 (i) in \cite{CPRY2019} that
	for any $k\geq 1$, the sequence
	$(\mathbb{W}_{k}^{(q), (n)})_{n \geq 1}$ and $(\mathbb{Z}_{k}^{(q), (n)})_{n \geq 1}$ of the corresponding scale functions in the bounded variation case (associated with $X^{(n)}$)  converges pointwise to scale functions $\mathbb{W}_{k}^{(q)}$ and $\mathbb{Z}_{k}^{(q)}$, respectively.
\end{proof}
\noindent Next, we shall use Theorem \ref{thm:refMultiResolvent} to get  an identity for the one-sided exit problem.
\begin{theorem}	\label{thm:OneSidedMultiRefRac}
	For any $q\geq 0$ and $x\leq a$, we have
		\begin{equation*}
			\Ex_{x}\bigl[e^{-q T_k^{a, +}}\mathbf{1}_{(T_k^{a, +} < \infty)}\bigr]=\frac{\mathbb{Z}_k^{(q)}(x)}{\mathbb{Z}_k^{(q)}(a)}.
		\end{equation*}
\end{theorem}
\begin{proof}
	Using the definition for $\mathbb{Z}_k^{(q)}(x)$ as in pp. 5445  in~\cite{CPRY2019}, we have
		\begin{equation}
			\mathbb{Z}_{k}^{(q)}(x):=1+q \sum_{i=0}^{k} \int_{b_i}^{b_{i+1}} \frac{\mathbb{W}_{k}^{(q)}(x ; y)}{\Xi_{\phi_i}(y)}\mathrm{d}y, \quad \text { where } \quad b_0=-\infty, \quad b_{k+1}=x,
			\label{eqn:Zinduction}
		\end{equation}
	or equivalently,
		\begin{align}
			\frac{\mathbb{Z}_{k}^{(q)}(x) -1}{q}
			=  \int_{0}^{\infty}\mathbb{W}_{k}^{(q)}(x ; y) \sum_{i=0}^{k}\frac{1}{\Xi_{\phi_i}(y)}\mathbf{1}_{(y \in (b_i, b_{i+1}])}\mathrm{d}y
			= \int_{0}^{x} \Xi_{\phi_k}(y)^{-1}\mathbb{W}_{k}^{(q)}(x ; y)\mathrm{d}y.
			\label{eqn:Zinduction2}
		\end{align}
	Then,  for  $e_q$ be an independent
	and exponentially distributed random variable with parameter $q >0$, it holds that
		\begin{align*}
			\Ex_{x}\bigl[e^{-q T_k^{a, +}}\mathbf{1}_{(T_k^{a, +} < \infty)}\bigr]
			= 1- q \int_{0}^{\infty}e^{-qt}\Prob_x(V_t \in \mathfrak{B}, t < T_k^{a, +} ) \mathrm{d}t
			=1-  q\Bigl[\frac{1}{q}\Bigl( 1- \frac{\mathbb{Z}_k^{(q)}(x)}{\mathbb{Z}_k^{(q)}(a)}\Bigr)\Bigr].
		\end{align*}
	where the last eq. follows from Theorem \ref{thm:refMultiResolvent} and use of Eq. \eqref{eqn:Zinduction2}.
\end{proof}
The rest of this section is devoted to  show that $V_k$ in Eq. \eqref{eqn:RefMulti} are monotonic based on their driving rate functions.
\begin{lemma}
	Suppose that for each $n \geq 1$, there exists a sequence a sequence of $(\phi_n)_{n\geq 1}$   with $\phi_n(x) \leq \phi_{n+1}(x)$ for all $x \in \R$. Then,
	$V_{n+1}(t) \leq V_{n}(t)$ for all $t \geq 0$.
	\label{lem:monotonocityVn}
\end{lemma}
\begin{proof} Amending  the arguments in \cite{CPRY2019} for the reflected case,  we consider a sequence of $(\phi_n)_{n\geq 1}$ satisfying the same conditions as in the proof of Proposition \ref{prop:convergenceVn} and
	we define the function \( \phi_{n+1}^{\varepsilon}(x):= \phi_{n+1}(x) + \varepsilon\)
	for  $\varepsilon > 0$. Then, we have $\phi_n(x) < \phi_{n+1}^{\varepsilon}(x)$
	for all $x \in \R$.
	Consider the process $V_{n+1}^{\varepsilon}(t)$ which is the solution to the following SDE
	\begin{align*}
			V_{n+1}^{\varepsilon}(t) &=
		 X(t) +  R^{\varepsilon}_{n+1}(t) - \int_{0}^{t}\phi_{n+1}^{\varepsilon}(V_{n+1}^{\varepsilon}(s))\mathrm{d}s,
		\end{align*}
	where $R^{\varepsilon}_{n+1}(t)$ has the same definition as in Eq. \eqref{eqn:RefMulti} for \(\phi_{n+1}^{\varepsilon}(x)\).
	
	Moreover, we define
	\begin{linenomath*}
		\[ \varsigma := \inf\{ t > 0: V_n(t) <  V_{n+1}^{\varepsilon}(t)\}, \]
	\end{linenomath*}
	and assume that $\varsigma < \infty$. Note that since $V_n$
	and $V_{n+1}^{\varepsilon}$ have the same jumps, the crossing of
	the path cannot occur at a jump instant.
	Thus, we obtain that \(V_n(t) - V_{n+1}^{\varepsilon}(t) \) is non-increasing in some
	$[\varsigma-\epsilon, \varsigma) $ for small enough $\epsilon$.
	Further, in $[\varsigma-\epsilon, \varsigma) $ we have also that \(R_n(t)-R_{n+1}^\varepsilon (t)\)
	  is increasing.
	Then, we have
		\begin{align*}
			0 \geq \frac{\mathrm{d}}{\mathrm{d}t}(V_n(t) -  V_{n+1}^{\varepsilon}(t))\Bigr|_{t = \varsigma^-}
			&=  \frac{\mathrm{d}}{\mathrm{d}t}\bigl(R_n(t)-R^\varepsilon_{n+1}(t)\bigr)\\
			&\quad + \phi_{n+1}^{\varepsilon}(V_{n+1}^{\varepsilon}({\varsigma^-)})) - { \phi_{n}(V_{n}(\varsigma^-))}> 0,
		\end{align*}
	which yields to
	$\varsigma = \infty$ which implies that $V_{n+1}^{\varepsilon}(t) \leq  V_{n}(t)$.
	Now, since $\phi_{n+1}(x) < \phi_{n+1}^{\varepsilon}(x)$ for all $x \in \R$, following the
	same argument as above, we obtain that $V_{n+1}^{\varepsilon}(t) \leq  V_{n+1}(t)$ for all $t \geq 0$.
	
	On the other hand, let
		\begin{align*}
			\Delta_{\varepsilon}(t) &:= V_{n+1}(t) - V_{n+1}^{\varepsilon}(t) \\
			&= R_{n+1}(t)-R^\varepsilon_{n+1}(t)+ \int_{0}^{t}\bigl[\phi_{n+1}^{\varepsilon}(V_{n+1}^{\varepsilon}(s))
			-\phi_{n+1}(V_{n+1}(s)) \bigr]\mathrm{d}s.
		\end{align*}
	Using integration by parts, noticing that \(	\Delta_{\varepsilon}(t)  \leq \int_{0}^{t}\bigl(\phi_{n+1}^{\varepsilon}(V_{n+1}^{\varepsilon}(s))
	-\phi_{n+1}(V_{n+1}(s)) \bigr)\mathrm{d}s\) and recalling
	that \(\phi_{n+1}^{\varepsilon}(x):= \phi_{n+1}(x) + \varepsilon\)
	gives that
		\begin{align*}
			(\Delta_{\varepsilon}(t))^2 &  \leq
			2 \int_0^t \Delta_{\varepsilon}(s)\bigl(\phi_{n+1}(V_{n+1}^{\varepsilon}(s))
			-\phi_{n+1}(V_{n+1}(s))\bigr)\mathrm{d}s
			+ 2\varepsilon\int_0^t \Delta_{\varepsilon}(s)\mathrm{d}s.
		\end{align*}
	Since $V_{n+1}^{\varepsilon}(t) \leq  V_{n+1}(t)$ for all $t \geq 0$ and $\phi_{n+1}$
	is a non-decreasing function,  we obtain that
		\[ \int_0^t \Delta_{\varepsilon}(s)\left(\phi_{n+1}\left(V_{n+1}^{\varepsilon}(s)\right)-\phi_{n+1}\left(V_{n+1}(s)\right)\right) d s \leq 0.\]
After these observations, we
get
	\[(\Delta_{\varepsilon}(t))^2 \leq 2 \varepsilon \int_0^t \Delta_{\varepsilon}(s)\mathrm{d}s \stackrel{\varepsilon \downarrow 0}{\longrightarrow} 0, \]
	which implies $\lim_{\varepsilon \downarrow 0}\Delta_{\varepsilon}(t) = \lim_{\varepsilon \downarrow 0}(V_{n+1}(t) - V_{n+1}^{\varepsilon}(t)) \equiv 0$. Hence,
		\(\lim_{\varepsilon \downarrow 0} V_{n+1}(t) = \lim_{\varepsilon \downarrow 0}V_{n+1}^{\varepsilon}(t) \leq \lim_{\varepsilon \downarrow 0} V_{n}(t),\)
	which completes the proof.
	\end{proof}
\subsection{Potential measure for reflected level-dependent L\'evy process}\label{A2}
In this section we provide the identity for  the resolvent of the (no-killed) reflected level-dependent L\'evy processes, defined in Eq. \eqref{eqn:SDELevelRef}. The following  result not only completes the literature of level-dependent L\'evy processes  (reflected at infima),  but also is used in Section \ref{killed reflection} to derive fluctuation identities for the  $\omega$-killed  reflected level-dependent L\'evy process.  The proof is similar to \cite{CPRY2019} adjusted for the reflected case and is provided below for completeness. Let $K_a$ have the same definition as in Section \ref{killed reflection}.
\begin{theorem}	\label{thm:resolventLevelRef}
	For $q\geq 0$, $x \leq a$, and a Borel set $\mathfrak{B}$ on $[0, a)$
	we have
		\begin{equation}			\label{eqn:resolventLevelRef}
			\Ex_x\Bigl[\int_0^{K_a} e^{-qt}\mathbf{1}_{(V_t \in \mathfrak{B})} \mathrm{d}t\Bigr]
			=\int_{\mathfrak{B}} \Xi_\phi(y)^{-1}\Bigl(\mathbb{W}^{(q)}(a; y)\frac{\mathbb{Z}^{(q)}(x)}{\mathbb{Z}^{(q)}(a)}-\mathbb{W}^{(q)}(x; y)\Bigr) \mathrm{d}y.
		\end{equation}
\end{theorem}
\begin{proof}
	Consider  the same approximating sequence $(\phi_{n})_{n\geq 1}$  as in the proof of Proposition  \ref{prop:convergenceVn} and $(V_n)_{n\geq 0}$ the corresponding sequence of non-increasing reflected
	multi-refracted L\'evy process (see also the proof of Proposition  \ref{prop:convergenceVn}). Also, let $\overline{V}_n(t) =  \sup_{0 \leq s \leq t}V_n(s)$ and $\overline{V}_t =  \sup_{0 \leq s \leq t}V_s$. For  $t \geq 0$, by triangle inequality we have  $\bigl|\overline{V}_n(t)-\overline{V}_t\bigr|  \leq \sup _{s \in[0, t]}\bigl|V_n(s)-V_s\bigr|$, along with  Proposition \ref{prop:convergenceVn},  yields that
	\begin{linenomath}
		\[\lim _{n \uparrow \infty}\left(V_n(t), \overline{V}_n(t)\right)=(V_t,  \overline{V}_t), \;\;\; \text{a.s.} \]
	\end{linenomath}
	Now, by using~Lemma \ref{lem:monotonocityVn}, for each $t >0$, the sequence
	$(V_n)_{n\geq 1}$ is non-increasing and thus for $a, y \geq 0$, 	$\{\overline{V}(t) \geq a\}=\bigcap_{n \geq 1}\{\overline{V}_n(t) \geq a\}$ and  $	\{V(t) \geq y\}=\bigcap_{n \geq 1}\left\{V_n(t) \geq y\right\}$.
	Therefore,  for any $x \in \R$, and $t > 0$,
		\begin{align*}
			\Prob_x(V_t \geq y, \overline{V}_t \geq a) & =\Prob_x\Bigl(\bigcap_{n \geq 1}\left\{V_n(t) \geq y, \overline{V}_n(t) \geq a\right\}\Bigr)  =\lim _{n \rightarrow \infty} \Prob_x\bigl(V_n(t) \geq y, \overline{V}_n(t) \geq a\bigr).
		\end{align*}
	Since  $\Prob_x(\overline{V}(t) = a) = 0$,  which follows verbatim from~\cite{PY2018},  implies that
		\begin{align*}
			\mathbb{P}_x(V_t \geq y, \overline{V}_t \leq a)
			&=  \mathbb{P}_x(V_t \geq y)
			-\mathbb{P}_x(V_t \geq y, \overline{V}_t \geq a) \\
			&=  \lim _{n \rightarrow \infty} \Prob_x\bigl(V_n(t) \geq y\bigr)
			-\lim _{n \rightarrow \infty} \mathbb{P}_x\bigl(V_n(t) \geq y, \overline{V}_n(t) \geq a\bigr) \\
			&= \lim _{n \rightarrow \infty}\Prob_x\bigl(V_n(t) \geq y, \overline{V}_n(t) \leq a\bigr).
		\end{align*}
	Thus, by Dominated Convergence Theorem we obtain that
		\begin{align*}
			\lim _{n \rightarrow \infty} \Ex_x\Bigl[\int_0^{T_n^{a,+}} e^{-q t} \mathbf{1}_{\left\{V_n(t) \geq y\right\}}\mathrm{d}t\Bigr]
			&=\lim _{n \rightarrow \infty} \int_0^{\infty} e^{-qt} \Prob_x\bigl(V_n(t) \geq y, \overline{V}_n(t) \leq a \bigr)\mathrm{d}t\\
			& =\int_0^{\infty} e^{-qt} \Prob_x\bigl(V_t \geq y, \overline{V}_t \leq a\bigr)\mathrm{d}t  =\Ex_x\Bigl[\int_0^{K_a} e^{-qt} \mathbf{1}_{\{V_t \geq y\}}\mathrm{d}t\Bigr] .
		\end{align*}
Recalling Eq.   \eqref{eqn:refMultiResolvent}, we have that
	\[	\lim _{n \rightarrow \infty} \Ex_x\Bigl[\int_0^{T_n^{a,+}} e^{-qt} \mathbf{1}_{\left\{V_n(t) \geq y\right\}}\mathrm{d}t\Bigr]= \lim _{n \rightarrow \infty} \int_{[y, \infty]} \Xi_{\phi_n}(z)^{-1}\Bigl(\mathbb{W}_n^{(q)}(a; z)\frac{\mathbb{Z}_n^{(q)}(x)}{\mathbb{Z}_n^{(q)}(a)}-\mathbb{W}_n^{(q)}(x; z)\Bigr)\mathrm{d}z.\]
	Then, by using the fact the scale functions $\mathbb{W}_n^{(q)}$ and $\mathbb{Z}_n^{(q)}$ converges to $\mathbb{W}^{(q)}$ and $\mathbb{Z}^{(q)}$, see Lemma 22 and Theorem 27 in \cite{CPRY2019},  and by the assumption of $\lim_{n\to \infty}\phi_{n} = \phi$ uniformly a.s., the result follows.
\end{proof}

\section*{Acknowledgement}
Z. Palmowski acknowledges that the research is partially supported by the Polish National Science Centre Grant No. 2021/41/B/HS4/00599. {M. \c{S}im\c{s}ek acknowledges  the financial support by the Scientific and Technological
Research Council of Turkey (T\"{U}B\.{I}TAK) through the B\.{I}DEB-2214/A
International Doctoral Research Fellowship Programme.}
The authors are grateful to the anonymous referees for their constructive comments and suggestions that have improved the content and presentation of this paper.

\bibliographystyle{plain}
\bibliography{Apostolos}
\end{document}